\documentclass[A4,11pt]{article}

\usepackage{color}

\usepackage{amsfonts}
\usepackage{amsmath}
\usepackage{amssymb}
\usepackage{algorithm2e}
\usepackage{graphicx}
\usepackage{float}
\usepackage[final]{pdfpages}
 \usepackage{hyperref} 
 \makeatletter
\hypersetup{
pdftitle={\@title},
pdfauthor={\@author},
colorlinks=true, 
breaklinks=true, 
urlcolor= blue, 
linkcolor= blue, 
citecolor=blue,    
bookmarksopen=true,
pdftoolbar=false,
pdfmenubar=false,
}
\makeatother

\newcommand{\eps}{\epsilon}

\newcommand{\qed}{$\square$}

\newcommand{\bfv}{{\bf v}}
\newcommand{\bfe}{{\bf e}}

\newcommand{\Q}{\mathcal{Q}}

\renewcommand{\Im}{\operatorname{Im} }
\renewcommand{\Re}{\operatorname{Re} }

\newcommand{\RR}{\mathbb{R}}
\newcommand{\CC}{\mathbb{C}}

\newcommand{\UU}{\mathbb{U}}

\newcommand{\LL}{\mathbb{L}}

\newcommand{\ux}{{\underline x}}

\newcommand{\caA}{{\cal A}}

\newcommand{\caD}{{\cal D}}
\newcommand{\caO}{{\cal O}}

\newcommand{\caE}{{\cal E}}
\newcommand{\caS}{{\cal S}}

\newcommand{\caZ}{{\cal Z}}

\newcommand{\caQ}{{\cal Q}}

\newcommand{\sfm}{{\mathsf m}}

\newcommand{\uOmega}{{\underline{{\caQ_\delta}}}}
\newcommand{\uu}{{\underline u}}

\newcommand{\uX}{{\underline X}}

\newcommand{\utau}{{\underline \tau}}

\newcommand{\ug}{{\underline g}}

\newcommand{\udelta}{{\underline \delta}}

\newcommand{\ou}{{\overline u}}

\newcommand{\JR}[1]{{\color{black} #1}}

\newtheorem{theorem}{Theorem}[section]

\newtheorem{remark}{Remark}[section]

\newtheorem{example}{Example}[section]

\newtheorem{proposition}[theorem]{Proposition}

\newtheorem{definition}[theorem]{Definition}

\newtheorem{corollary}[theorem]{Corollary}

\newtheorem{lemma}[theorem]{Lemma}

\numberwithin{equation}{section}

\begin{document}

\title{
Perfectly Matched Layers on Cubic Domains for 
Pauli's Equations
 }

\date{}
\author{ Laurence Halpern
\footnote{Universit\'e Sorbonne Paris Nord, LAGA, CNRS, UMR 7539,  F-93430, Villetaneuse, France,
halpern@math.univ-paris13.fr}
\and 
{Jeffrey  Rauch} 
\footnote{Department of Mathematics,
University of Michigan, Ann Arbor  48109  MI, USA
rauch@umich.edu.}
}

\maketitle

\begin{abstract} 

This article proves the wellposedness of 
the boundary value problem that arises when    
PML algorithms are applied 
to Pauli's equations
with a three dimensional rectangle as computational domain.
The absorptions are positive near the 
boundary and zero far from the boundary so are always 
$x$-dependent.
At the flat parts of the boundary of the rectangle,  the natural absorbing boundary conditions 
are imposed.  The difficulty addressed is the analysis of the resulting
variable coefficient problem on the rectangular solid with its edges and corners.
The Laplace transform is analysed.  
\JR{We derive an additional
boundary condition that is automatically satisfied and yields 
a coercive Helmholtz boundary value problem  on smoothed
boundaries with uniform estimates justifying the limit
of vanishing smoothing.
}
 This yields
  the first stability proof 
with  $x$-dependent absorptions 
on a domain whose boundary is not smooth.
\end{abstract}



{\bf Acknowledgements.}  \JR{We thank the  referee for suggestions
that greatly improved the presentation. }   JBR gratefully acknowledges the support of the 
CRM Centro De Giorgi in Pisa, and the LAGA  at 
the Universit\'e Sorbonne Paris Nord for support of this research.

 \break
 \tableofcontents

  \section{Introduction}

  This paper analyses initial boundary value problems that arise
  when one uses perfectly matched absorbing layers
  in the time domain.
    The most common configuration is a three dimensional rectangular solid 
  surrounded  by a larger rectangular  solid computational domain.  The inner solid
  contains the sources and is the region where the computed values are 
  required.
  In the region between the rectangles,  perfectly
  matched layers are interposed.
  Boundary conditions at the exterior boundary
  are imposed that are designed to be weakly reflecting.  
  In addition to perfect matching, an advantage of the PML
  strategy is its ease of implementation including at the corners.
  To our knowledge, the present work is the first to prove
  wellposedness for such a PML  with non constant 
  absorptions
  $\sigma_j$ in the presence of trihedral corners.  That problem
  poses two fundamental
    challenges.
      
  Even for a system with a very simple energy estimate
  like Pauli's equations, the split equations of B\'erenger and also the 
  stretched  system that is at the heart of its analysis do not 
  have simple estimates.       Such estimates are crucial for 
  constructing solutions and express stability.   In practice the 
  split system needs to be discretized and the stability of the 
  discretization analysed. This article does not study that 
  problem.  A recent survey for the constant coefficient
  half space case is \cite{coulombel2019transparent}.

 The Pauli  system shares the Lorentz invariance, symmetry, and  three dimensionality
 of Maxwell's equations.   It has two advantages.   It is a $2\times 2$ system
 as opposed to a $6\times 6$ system.  More importantly,
 the generator is elliptic.  The analysis
 extends with almost no modifications to the Dirac system.
 The Maxwell system poses serious problems.  It's treatment is work in 
 progress.
 The Pauli  operator is
  \begin{equation}
  \label{eq:paulisystem1}
 L := \partial_t
+
\begin{pmatrix}
1 & 0
\cr
 0 & -1
 \end{pmatrix}
 \partial_1
  +
 \begin{pmatrix}
0 & 1
\cr
 1 & 0
 \end{pmatrix}
 \partial_2
 +
  \begin{pmatrix}
0 & i
\cr
 -i & 0
 \end{pmatrix}
 \partial_3
  :=
 \partial_t + \sum_{j=1}^3
 A_j\partial_j.
 \end{equation}
 Introduce the notations with $\xi\in \CC^3$,
 \begin{equation}
 \label{eq:paulisystem2}
 L(\partial_t,\partial_x)
 \ :=\ 
 \partial_t \ +\ A(\partial_x),
 \qquad
 A(\xi) \, :=\, A_1\xi_1 \, +\, A_2\xi_2\, +\,  A_3\xi_3\,.
  \end{equation}

\begin{definition}
\label{def:caE}
For $A\in {\rm Hom}(\CC^k)$, with  spectrum disjoint from $i\RR$,
$\caE^+(A)$
(resp.     $\caE^-(A)$)    denotes the spectral subspace corresponding to eigenvalues with 
strictly
positive
(resp.   strictly negative) real part.
Denote by  $\pi^\pm(A)$  the corresponding spectral
projections onto those spaces.
As our interest is the  Pauli system,  $\caE^\pm(\xi)$
and
$\pi^\pm(\xi)$
are
 shorthands
for $\caE^\pm(A(\xi))$ and
$\pi^\pm(A(\xi))$ for $\xi\in \CC^3$ so that   $A(\xi)$ 
has no purely imaginary eigenvalues.
\end{definition}

\begin{definition}
Denote by $\caQ=\caQ(L_1,L_2,L_3)$ the  rectangle 
$$
\caQ \ :=\ \big\{x\in \RR^3 \ :\  |x_j|<L_j/2,\ \ j=1,2,3
\big\}  .
$$
$\caQ$ has six open faces $G_k$ with $1\le k\le 6$.  For $j=1,2,3$,
\begin{align*}
&G_j:= 
\big\{x_j=-L_j /2,\ 
{\rm and},\ 
 |x_i|<L_i /2 \ {\rm for}\  j\ne i
 \big\},
 \cr
  &G_{j+3}:= 
\big\{x_j=L_j/2,\ 
{\rm and},\ 
 |x_i|<L_i/2 \ {\rm for}\   j\ne i
 \big\}.
\end{align*}
For a point  $x\in G_k$,
 $\nu(x)$ denotes the outward unit normal to $\caQ$.
The split equations involve non negative 
{\bf absorption coefficients} $\sigma_j\in C^\infty_0(\RR)$ for $j=1,2,3$.
\end{definition}

\begin{example}  
\label{ex:standard}
For the {\bf usual implementations of the PML method}, 
there is an $\ell<1$ so that the absorptions vanish in $\ell\Q$,
the sources are supported in $\ell\Q$ and the values
of the solution on $\ell\Q$ are those of interest.
\end{example}

The Pauli system is 
  \begin{equation}
  \label{eq:Pauli10}
 \Big(\partial_t
\ +\
A_1 \partial_1
\  +\
A_2
 \partial_2
 \ +\
  A_3
 \partial_3
 \Big)
 u
 \ =\ f
 \quad
 {\rm on}\ 
 \ \caQ.
 \end{equation}
 
\begin{definition}
B\'erenger's method  has unknown that is 
 a triple 
$ (U^1,U^2, U^3)$ with $U^j$ taking
values in $\CC^2$.
 On
 $\RR\times \caQ$,
 $(U^1,U^2, U^3)$ satisfy 
  the {\bf split equations},
 \begin{equation}
\label{eq:splitpauli}
\begin{aligned}
\big(\partial_t + \sigma_1(x_1)\big)U^1 \ +\ A_1\partial_1\Big(U^1+U^2+U^3\Big)\ &=\ f_1\,,
\cr
\big(\partial_t + \sigma_2(x_2)\big) U^2 \ +\ A_2\partial_2\Big(U^1+U^2+U^3\Big)\ &=\ f_2\,,
\cr
\big(\partial_t + \sigma_3(x_3)\big) U^3 \ +\ A_3\partial_3\Big(U^1+U^2+U^3\Big)\ &=\ f_3\,.
\end{aligned}
\end{equation}
The $j^{\rm th}$ equation has the $\partial_j$ derivative.
The $f_j$ are constrained to satisfy $f=\sum_j f_j$ and 
to vanish on a neighborhood of $\partial\caQ$.  A
choice respecting the symmetry of the problem is
 $f_j=f/3$ for $j=1,2,3$.
 \end{definition}

  The boundary of $\caQ$ is not perfectly transparent.  
  In favorable cases like the Pauli
  system, waves are expected to decay  in the layers
  so  little signal
  reaches $\partial\caQ$ and the reflections cause small errors.  In
  practice rather thin layers suffice.   
  With $x$-dependent absorptions and computations in the time domain,
  proving exponential decay in the layers is an outstanding open
  problem.
  See Remark \ref{rem:timeharmonic} for the easier  time harmonic case.

The  split equations  are not
symmetric and they  have a lower order term that depends on $x$ through
the absorption coefficients $\sigma_j$. They do not have simple {\it a priori} estimates
showing that they yield a well posed pure initial value problem.
Petit-Bergez 
 \cite{Petit:2006:PFB, halpern2011analysis} 
 proved that since the Pauli system generates a $C_0$-semigroup
 on $L^2(\RR^3)$
 and has  elliptic generator, it follows 
  that the split equations on $\RR^3$ also generate 
a $C_0$-semigroup.   This contrasts to the
loss of one derivative for the split Maxwell equations proved by
Arbarbanel and Gottlieb
\cite{Abarbanel:1998:WPM}.

The Pauli system is symmetric hyperbolic.   The most 
strongly dissipative boundary condition for the Pauli system  is $u\in \caE^+(\nu)$.
Thanks to the symmetry and ellipticity
 there is an $M_0$ so that for   $M>M_0$ and $f\in e^{Mt} L^2(\RR\times \caQ)$
the boundary value problem $Lu=f$ with boundary condition
$u\in \caE^+(\nu)$ on the $\RR\times G_k$
 has a unique
solution $u\in e^{M t}L^2(\RR\times\caQ)$
with $u\in e^{Mt}L^2(\RR\times \partial\caQ)$
(see Part I of \cite{HR:2016:HBV}).  The substantial
difficulty here is the presence of the edges and corner
of $\Q$.

For the Pauli system, if one imposed
a conservative rather than dissipative boundary  condition,
then waves that arrive at the external boundary 
would be totally reflected
back   to the interior.   Approximately this  behavior would be inherited
by the split system and is to be avoided.   
Less obvious
is that our proof breaks down for non dissipative conditions.
Even for the unstretched problem uniqueness of solutions
for a conservative condition on domains with trihedral corners
is not known.
This underscores the  difficulties  of 
domains with
trihedral corners.

The Laplace transform of solutions of 
\eqref{eq:splitpauli} satisfy
  \begin{equation}
\label{eq:split10}
\begin{aligned}
\big(\tau+ \sigma_1(x_1)\big) \widehat U^1 \ +\ A_1\partial_1\Big(\widehat U^1+
\widehat U^2+\widehat U^3
\Big)\ &=\ \widehat f_1\,,
\cr
\big(\tau + \sigma_2(x_2)\big)\widehat U^2 \ +\ 
A_2\partial_2\Big(\widehat U^1+ \widehat  U^2+\widehat U^3\Big)\ &=\ \widehat f_2\,,
\cr
\big(\tau + \sigma_3(x_3)\big)\widehat U^3 \ +\ 
A_3\partial_3\Big(\widehat U^1+ \widehat  U^2+\widehat U^3\Big)\ &=\ \widehat f_3\,.
\end{aligned}
\end{equation}

\JR{
\begin{definition} 
Define for $j=1,2,3$,
\begin{equation}\label{eq:deftilde}
\widetilde\partial_j
\ :=\
\frac{\tau}{\tau + \sigma_j(x_j)}
\,
\frac{\partial}
{\partial x_j}
\,,\qquad {\rm and } \quad\
u\ :=\ 
 \widehat U^1+ \widehat  U^2+ \widehat  U^3\,.
\end{equation}

The  {\bf stretched Pauli operator} is defined by
$$
L(\tau, \widetilde \partial_x)\ :=\
 \tau 
\ +\ 
\sum_{j=1}^3
A_j\widetilde\partial_j
\ =\ 
 \tau 
\ +\ 
\sum_{j=1}^3
A_j
\frac{\tau}{\tau + \sigma_j(x_j)}
\,
\frac{\partial}
{\partial x_j}\,.
$$

The stretched operator yields 
 {\bf the stretched equation}
\begin{equation}
\label{eq:bermagic10}
\Big(\tau 
\ +\ 
A_1\widetilde\partial_1
\ +\ 
A_2\widetilde\partial_2
\ +\ 
A_3\widetilde\partial_3
\Big)
u \ =\ 
\ F\ :=\
\sum_{j=1}^3 \frac{\tau}{\tau + \sigma_j(x_j)} \hat{f}_j \,.
\end{equation}
\end{definition}

These objects are called stretched because when $\tau$
is real they arise from a change of variable in $x$, 
see Section \ref{sec:stretch}.

}

\begin{definition}
{\bf i.}  If $\caO\subset\RR^3$ is open and $K\subset \caO$ is compact,
$$
C^\infty_K(\caO)
\ :=\
\Big\{ f\in C^\infty(\caO)\ ;\ 
{\rm supp}\,f\subset K
\Big\}.
$$
{\bf ii.}  Similarly,
$$
L^2_K(\caO)
\ :=\ 
\Big\{
f\in L^2(\caO)\ ;\ 
{\rm supp}\,f\subset K
\Big\}\,.
$$
\end{definition}

The next result solving the stretched equation is the main result of the paper.
It allows one to prove the stability of B\'erenger's split method.

\begin{theorem}
\label{thm:mainstretched}
For each $\ell\in ]0,1[$ 
there exist $C,M_1$ 
so that for all $M\ge M_1$, and,  holomorphic
$F:{\{\rm Re}\, \tau>M\} \to  L_{\ell\overline{\Q}}^2(\caQ)$, 
there is a unique holomorphic function $u:\{{\rm Re}\,\tau>M\}\to H^1(\caQ)$
 satisfying the stretched boundary  value problem
  on $\caQ$,
 \begin{equation*}
\label{eq:stretchedpair2}
L(\tau,\widetilde\partial_x){u}\ =\ F
\quad  {\rm on} \quad \caQ,
\qquad
{u} \, \in \, \caE^+\big(\nu \big)
\quad
{\rm on}
\quad
G_k, \ \  1\le k\le 6.
\end{equation*}
It  satisfies for all ${\rm Re}\,\tau>M$,
\begin{align}
\label{eq:arnold}
({\rm Re}\, \tau)\,
\big\|{u}
\big\|_{L^2(\caQ)}
+
({\rm Re}\, \tau)^{1/2}\,
\big\|
 u
 \big\|_{L^2(\partial\caQ)}
 +
 \frac{\Re\tau}{|\tau|}\,
 \| \nabla_xu\|
\,\le\, 
C\, \big\|
F
\big\|_{L^2(\caQ)}
.
\end{align}
\end{theorem}
\begin{remark}
{\bf i.}
The gradient estimate degenerates as $\Im\tau \to \infty$ with $\Re\tau$ fixed.
A second hyperbolic aspect of \eqref{eq:arnold}  is that the boundary values are estimated in 
$L^2$ and not in $H^{1/2}$.

\JR{
{\bf ii.}   The holomorphy is crucial.   The theorem is used to 
construct the Laplace transform of an object supported in $t\ge 0$ 
that
must be holomorphic.    In addition, uniqueness 
is reduced by an analyticity argument
to our uniqueness theorem for symmetric hyperbolic
problems in domains with trihedral corners  \cite{HR:2016:HBV}.
}
\end{remark}

 \begin{theorem} 
\label{thm:pauli}   There are strictly positive constants $C,M$
so that if $\lambda>M$,
 ${f}\in e^{\lambda t}L^2(\RR\times\caQ) $ 
with support in 
$[0,\infty[\times \ell\overline\Q$, then there is one and
only one  
$
\big(U^1,U^2,U^3)\in e^{\lambda t}L^2(\RR\times\caQ)
$
 supported in $t\ge 0$
that satisfies
 \eqref{eq:splitpauli},  
 and the boundary condition $U^1+U^2+U^3\in \caE^+(\nu)$
 on  each $G_k$.
The function $U^1+U^2+U^3$ 
satisfies,
\begin{equation}
\label{eq:berpauliest1}
\begin{aligned}
\lambda\,   \big\| & e^{-\lambda t}  
(U^1+U^2+U^3)
\big\|_{L^2(\RR\times  \caQ)} \ +
\cr
&
\lambda^{1/2}\big\| e^{-\lambda t} 
(U^1+U^2+U^3)
\big\|_{L^2(\RR\times\partial\caQ)}
\, \le\,
  C\,
  \big\|  e^{-\lambda t}   \,f 
\big\|_{L^2(\RR\times\caQ)}.
\end{aligned}
\end{equation}
The split unknowns satisfy the weaker estimate
\begin{equation}
\label{eq:berest1}
\begin{aligned}
\big\| e^{-\lambda t} \{\lambda U^j\,,\,\partial_{t}U^j\}
\big\|_{L^2(\RR: H^{-1}(\caQ))}
\ \le\
C\,
\big\| e^{-\lambda t} \, f
\big\|_{L^2(\RR\times \caQ)}
\,.
\end{aligned}
\end{equation}
\end{theorem}

Theorem \ref{thm:pauli}
allows us to analyse the 
split equations with the absorbing boundary  condition 
$U^1+U^2+U^3\in  \caE^+(\nu)$ on the $G_k$.
 It is the first existence theorem for the split equations with non constant
$\sigma_j$ in domains whose boundary is not smooth.   Since 
standard practice uses 
cubes with
non constant $\sigma_j$ is it is the first justification, beyond
extensive practical experience,   
 that the  B\'erenger algorithm is stable.

\begin{remark}
{\bf i.}  
It is wise to think of B\'erenger's algorithm as
a method that inputs 
$f$ and outputs $U^1+U^2+U^3$.   Estimate
\ref{eq:berpauliest1}
shows that the output  satisfies 
 bounds as strong as strictly dissipative boundary value
 problems for symmetric hyperbolic systems.   This behavior
  is known for the pure initial value problem 
 (see  Theorem 1.3 of 
 B\'ecache and Joly 
 \cite{Becache:2002:ABP}  for 
 the split  Maxwell equations with constant $\sigma_j$,
 and 
 \cite{halpern2011analysis}
 for B\'erenger transmission problems  with 
 variable $\sigma_j$).

{\bf ii.}   The estimates of Theorem \ref{thm:pauli}  
permit exponential growth in time. 
Even for sources compactly supported in time.    Practical experience
with B\'erenger's method for equations 
closely tied to the wave equation (e.g. Maxwell and Pauli) show no growth in time
even with variable $\sigma_j$.   Interesting bounds uniform in time
are proved for the case of  constant $\sigma_j$ 
  for 
sufficiently regular  
solutions by 
B\'ecache-Joly,  Diaz-Joly, and Baffet-Grote-Imperiale-Kachanovska
\cite{Becache:2002:ABP, Diaz:2006:TDA,baffet2019energy}.    
  Uniform bounds in time 
 is an important and wide open 
problem for variable $\sigma$ even for the problem on $\RR^{1+d}$.
Appelo-Hagstrom-Kreiss  \cite{Appelo:2006:PML}
analyse the problem of exponential growth with constant parameters by explicit
formulas in  Fourier.  They propose stabilization
methods.    Variable coefficients and corner domains 
are beyond
 that strategy.
 
 {\bf iii.}  \JR{Other versions of PML lead to the same stretched system.
 The stability theorem for the stretched boundary value problem
implies stability for these versions too.
Respecting the history, we present the details
  for   B\'erenger's splitting.}
\end{remark}

The paper is organized as follows.  Section \ref{sec:Pauli} 
presents the Pauli system and most importantly the stretched
Pauli system that is satisfied by the Laplace transform of
$\widehat U^1 +\widehat U^2 +\widehat U^3$.  
\JR{
Theorems
\ref{thm:mainstretched} 
and
\ref{thm:deltastretched} 
assert existence and uniqueness for the boundary value problems
for the stretched system on $\caQ$ as well as smoothed versions
on
$\caQ_\delta$.  It is crucial that these results are proved 
with  $\delta$-independent estimates that justify passing to the limit
$\delta\to 0$.   
}

\JR{
 It is routine to show that solutions of the stretched
system are solutions of a Helmholtz type equation.   
An  important step is showing that the solutions
on $\caQ_\delta$ satisfy an additional boundary condition
stated in Corollary 
\ref{cor:BC2}.

The second boundary condition yields a 
coercive elliptic boundary value problem that is studied in 
Section \ref{sec:PauliHelmholtz}.
 Theorem \ref{thm:perturbed} yields the important
 uniform estimates for this boundary value problem.
 They are derived 
  by the energy method tied
to a  family of complex quadratic forms.   The
real and imaginary parts play key roles. The 
geometry  is singular in 
the limit $\delta\to 0$.  In spite of this, $H^1$ estimates uniform 
in $\delta$ and $\tau$ are proved.  The $H^2$ estimates degenerate when $\delta\to 0$.}

We have considered the option of skipping the smoothing
and using layer potential methods 
developed for the study of Lipschitz domains.   
Since the hard harmonic analysis would need to be 
adapted to the new problems,
the smoothing is both more elementary and shorter.

   \JR{
Section \ref{sec:maintheorems}  derives the main theorems
from the Helmholtz existence results. 
  Section \ref{sec:stretcheddelta}
proves unique solvability of the 
stretched Pauli system on $\Q_\delta$
stated in 
Theorem \ref{thm:deltastretched}.

First it is proved that the solution of the stretched Helmholtz
boundary value problem on $\Q_\delta$ satisfies the  stretched Pauli 
boundary value problem.
Here the $H^2$ smoothness of solutions 
on $\caQ_\delta$
is important.  Then the $H^1$ limit $\delta\to 0$
yields solutions of the  stretched Pauli system on $\caQ$.
 Holomorphy in $\tau$ is crucial
for uniqueness.

Section
\ref{sec:stretched}  proves Theorem \ref{thm:mainstretched}
asserting solvability of the stretched Pauli boundary  value problem on
$\caQ$ 
by passing to the limit $\delta \to 0$. 
 Section
\ref{sec:splitpauli} derives 
Theorem \ref{thm:pauli} asserting the solvability of the split equations
by constructing the Laplace transform by solving a stretched Pauli system.
}

The  proof is long and technical.  The hypothesis $\sigma_j\in C^\infty$ 
  avoids some inessential difficulties.   The proof  uses
  $H^2$ regularity for the Helmholtz problem on
  $\caQ_\delta$.  Absorptions $\sigma_j\in L^\infty$ suffice for 
  $H^1$ regularity and $\sigma_j$  lipschitzian is sufficient
  for $H^2(\Q_\delta)$.
  Standard practice involves such lipschitzian absorptions. 
  This  strengthening of the results is left to the interested reader.

\section{The Pauli system and smoothed domains $\Q_\delta$}
\label{sec:Pauli}

 \subsection{Pauli system and its  symbol}
 \label{sec:paulieig}
 
 The coefficients of the Pauli system \eqref{eq:paulisystem1} satisfy,
  \begin{equation}
\label{eq:anticommutation}
A_j^2\ =\ I,
\qquad
A_iA_j + A_jA_i \ =\ 0 \quad {\rm for}
\quad
i\ne j\,.
\end{equation}
These identities imply the connections to the Laplacian,
\begin{equation}
\label{eq:anticomm1}
\Big(
\sum_j A_j\partial_j\Big)^2
\ =\ \Delta \,,
\quad\ \
\Big(\sum A_j\partial_j-\tau\Big)
\Big(\sum A_j\partial_j+\tau\Big)
\ =\ \Delta -\tau^2\,.
\end{equation}

\begin{proposition}
\label{prop:paulispec}
\JR{With  $L$ and $A$ from the Pauli operator 
\eqref{eq:paulisystem2}, and the conventions  of
Definition \ref{def:caE}, the following hold.}

{\bf i.}  For all $(\tau,\xi)\in  \CC^{1+3}$,
  \begin{equation}\label{eq:detL}
  \det L(\tau,\xi) \ =\ \tau^2\ -\ 
  \sum_{j=1}^3 \xi_j^2\,.
  \end{equation}

{\bf ii.}  For $\xi\in\RR^3\setminus 0$, the $2\times 2$ hermitian symmetric matrix $A(\xi)$ has 
  eigenvalues $\pm |\xi|$ with one dimensional eigenspaces 
$$
\caE^-(\xi  ) \, =\,
\CC
\big(\xi_1- | \xi|\,,\,
\xi_2-i\xi_3
\big),
\quad
{\rm and},
\quad
\caE^+(\xi  ) \, =\,
\CC
\big(\xi_1+ | \xi|\,,\,
\xi_2-i\xi_3
\big).
$$

{\bf iii.}   For  all $\xi, \eta\in \CC^3$,
 \begin{equation}
 \label{eq:anticomm2}
 A(\xi)A(\eta)+A(\eta)A(\xi)
\ =\ 
2\Big(\sum_i  \xi_i\eta_i\Big)I\,.
\end{equation}
\end{proposition}

{\bf Proof.}  Write from \eqref{eq:paulisystem2},
$$
L(\tau,\xi)=\tau +A(\xi)=
\begin{pmatrix}
\tau+\xi_1 & \xi_2+i\xi_3\\
\xi_2-i\xi_3&\tau-\xi_1 
\end{pmatrix}
$$
This implies 
 {\bf i}.

{\bf ii.}  
For $\xi\in \RR^3\setminus 0$, \eqref{eq:detL} shows that the eigenvalues of $A(\xi)$ are  $\pm |\xi|$. The   first column yields the formula for $\caE^-(\xi)$ in {\bf ii.}
The other choice of sign yields $\caE^+(\xi).$

{\bf iii.}  Expand
$$
A(\xi)A(\eta)=\Big(
\sum_iA_i\xi_i
\Big)
\Big(
\sum_jA_j\eta_j
\Big)
\ =\ 
\sum_{i,j}  A_iA_j \xi_i\eta_j\,.
$$
Symmetrizing yields,
$$
A(\xi)A(\eta)+ A(\eta)A(\xi)\ =\
\sum_{i,j}  A_iA_j \xi_i\eta_j
\ +\
\sum_{i,j}  A_iA_j \eta_i\xi_j\,.
$$
In the last sum interchange the role of $i,j$  to find
$$
A(\xi)A(\eta)+ A(\eta)A(\xi)\ =\
\sum_{i,j}  A_iA_j \xi_i\eta_j
\ +\
\sum_{i,j}  A_jA_i \eta_j\xi_i\,.
$$
Separate out the terms with $i=j$ to find
$$
A(\xi)A(\eta)+ A(\eta)A(\xi)
\ =\ 
2\sum_i A_i^2\,\xi_i\eta_i
\ +\
\sum_{i\ne j}
\Big(
A_iA_j 
\ +\
 A_jA_i 
 \Big)
 \eta_i\xi_j\,.
$$
 Equation \eqref{eq:anticommutation} 
 yields
 \eqref{eq:anticomm2}.
\hfill
\qed
\vskip.2cm

\begin{example}
\label{ex:lambda}
  Define
$
\caZ :=
\big\{\xi\in \CC^3 : \sum_j \xi_j^2 =0\big\}$.
For $\xi\in \CC^3\setminus \caZ$, 
the spectrum of $A(\xi)$ consists of two simple eigenvalues differing
by a factor $-1$.  

The eigenvalues $\pm |\xi|$ for $\xi\in \RR^3\setminus 0$
extend  to holomorphic
eigenvalues $\lambda^\pm(\xi)=\pm (\sum \xi_j^2)^{1/2}$
on  the domain
\begin{equation}
\label{eq:holdomain}
\Big\{
\xi\in \CC^3\setminus 0\,:\,
|{\rm Im}\,\xi|<|{\rm Re}\,\xi|
\Big\}
.
\end{equation}
In this case $\sum \xi_j^2$ belongs to the simply connected
subset $\CC\setminus ]-\infty,0]\subset \CC\setminus 0$.
\end{example}

 \begin{proposition}  
 \label{prop:paulispec2}
  {\bf i.}  The 
eigenprojections $\pi^{\pm}(\xi)$  for 
 $\xi\in \RR^3\setminus 0$ extend  to holomorphic functions
on  the domain 
\eqref{eq:holdomain}, 
satisfying with notation from Example \ref{ex:lambda},
\begin{equation}
\label{eq:complexev}
\pi^\pm(\xi)\,A(\xi) 
\ =\
A(\xi) \,\pi^\pm(\xi)
\ =\
\lambda^\pm(\xi)\ 
\pi^\pm(\xi)\,.
\end{equation}

 They are given by
 $$
 \pi^\pm(\xi)
 \ =\
 \frac12\
  \Big(\sum \xi_j^2\Big)^{-1/2}
  \Big(
 A(\xi) \pm 
\big(\sum \xi_j^2\big)^{1/2}
 \,I 
 \Big)
 \,.
 $$
 
{\bf ii.}   For  $\xi, \eta$ belonging to \eqref{eq:holdomain},
 \begin{equation}
 \label{eq:piApi}
\pi^\pm(\eta)\, A(\xi)\, \pi^\pm(\eta)\ =\ 
  \big(\sum \eta_j^2\big)^{-1/2}
\,
\big(\sum \xi_i\eta_i\big)\,
\pi^\pm(\eta)\,.
\end{equation}
\end{proposition}

{\bf Proof.}   {\bf i.}
The formulas
$$
A(\xi) =\lambda^+ \big(\pi^+(A(\xi)) -  \,\pi^-(A(\xi))\big),
\quad
{\rm and},
\quad
I \ =\ 
\pi^+(A(\xi)) + \pi^-(A(\xi)),
$$
together with
$(\pi^\pm(\eta))^2=\pi^\pm(\eta)$, imply the formulas for $\pi^\pm(A(\xi))$ in {\bf i.}

{\bf ii.}
Multiply \eqref{eq:anticomm2} on the left and right by 
$\pi^\pm(\eta)$ 
to find
$$
2\,\pi^\pm(\eta)\Big(\sum \xi_i\eta_i\Big)\pi^\pm(\eta)
\, =\,
\pi^\pm(\eta)A(\xi)A(\eta)\pi^\pm(\eta)
\, +\,
\pi^\pm(\eta)A(\eta)A(\xi)\pi^\pm(\eta)
\,.
$$
Use \eqref{eq:complexev} twice and
$(\pi^\pm(\eta))^2=\pi^\pm(\eta)$ to find,
\begin{align*}
2\,\Big(\sum \xi_i\eta_i\Big)\pi^\pm(\eta)
\, &=\,
\pi^\pm(\eta)A(\xi)\lambda^\pm(\eta)\,\pi^\pm(\eta)
\ +\ 
\lambda^\pm(\eta)\,\pi^\pm(\eta)A(\xi)\pi^\pm(\eta)
\cr
\  &=\
 2\,\lambda^\pm(\eta)\,
\pi^\pm(\eta)A(\xi)\pi^\pm(\eta)
\,.
\end{align*}
This completes the proof.
 \hfill
 \qed
 \vskip.2cm

\subsection{\JR{ The stretched system on smoothed domains $\caQ_\delta$}}
\label{sec:stretch}

The stretched equation \eqref{eq:bermagic10} resembles the Laplace
transform of the  original system. For $\tau$ real and positive
it comes from the original
transformed system by a change of variable, called 
{\it coordinate stretching} (see Section \ref{sec:stretchedPauli},
and  Chew-Weedon \cite{Chew:1994:IIA}).

\JR{
 \begin{definition}
 \label{def:Xj}
{\bf i.}    For $\tau\in \CC\setminus \{0\}$ the 
 coordinate stretchings  
 $X_j(\tau,x_j)$
 are defined  as  the solutions of the ordinary differential equation in $x_j$,
  \begin{equation}
  \label{eq:defXj}
  \frac{\partial X_j}{\partial  x_j} \ =\ 
  \frac{\tau + \sigma_j(x_j)}
  {\tau}\,,
  \qquad
  X_j(0)\,=\,0\,.
  \end{equation}
    
  {\bf ii.}  For real $\tau>0$,
  $\partial_jX_j>0$ and    $x\mapsto X(\tau,x)$ is a diffeomorphism from $\RR^3$
  onto itself.
     Denote by 
  $\uOmega\subset\RR^d_X$ 
   the image of 
    ${Q_\delta}\subset \RR^d_x$.  
  \end{definition}
  
  \begin{example}
  In the standard implementation 
  of Example \ref{ex:standard},  the 
  $\sigma_j$ vanish on 
  $\ell\Q$.
  Therefore 
  $X$
is equal to the identity on that set.
  \end{example}

  Compute for real $\tau>0$,
\begin{equation}
\label{eq:cov}
\frac{\partial}{\partial x_j}   =\sum_k
\frac{\partial X_k}{\partial x_j}\,
\frac{\partial}{\partial X_k}
 =
\frac{\tau + \sigma_j(x_j)}{\tau}\,
\frac{\partial}{\partial X_j},
\quad
  \frac{\tau}{\tau+ \sigma_j(x_j)} \, \frac{\partial} {\partial  x_{j}}
   =
  \frac{\partial}{\partial X_j}
  .
\end{equation}
   Equation \eqref{eq:cov} gives a geometric interpretation of the stretched
   operator $L(\tau,\widetilde\partial)$ for $\tau\in \RR_+$.   It shows that $\widetilde \partial_j$ in the 
   $x$ coordinates is equal to $\partial/\partial X_j$ in the $X$ coordinates.
   Therefore if $u(x)$ and $v(X)$ are related by $v(X(\tau,x)) =u(x)$ then
   $
   L(\tau, \widetilde \partial) u(x) =\big(L(\tau, \partial_X) v\big)(X(\tau,x))$.
   
}

The stretched equations are sometimes expressed using auxiliary
variables $\psi_j$ defined as the solutions of 
$$
(\partial_t+\sigma_j(x_j) )\psi_j
\ =\ \partial_t u,
\qquad
\psi_j =0\quad
{\rm for}
\quad t<0.
$$
Then
$
\widetilde \partial_j \widehat u \  =\ 
\partial_j \widehat \psi_j\,.
$

Theorem 
\ref{thm:mainstretched} 
is proved by solving the stretched equation
on smoothed truncated domains and passing to the limit.

 \begin{definition}
 \label{def:Odelta}
 The singular set of the boundary of $\cal Q$ is 
 $$
 \caS \ :=\
 \Big\{
 x\in \partial {\cal Q}\,:\,  \exists i\ne j, \ x\in \overline G_i \cap \overline G_j
 \Big\}.
 $$

 Introduce for $0<\delta<1$ bounded smooth  approximations $\caQ_\delta$
 of $\caQ$.  
  Smooth the edges and corners of $\caQ$
 on a $\delta/2$-neighborhood of $\caS$ to yield bounded smooth convex 
 sets $\caQ_\delta$.  Do this so that for $\delta_1<\delta_2$,
 $\caQ_{\delta_1}\supset \caQ_{\delta_2}$.
  \end{definition}
 
  \begin{definition}
For $\tau$ with  ${\rm Re}\,\tau >0$ and $\nu\in \RR^3$  define
$$
\widetilde \nu (\tau,x)
\ :=\
\Big(
\frac{\nu_1\, \tau}{\tau+ \sigma_1(x)}
\,,\,
\frac{\nu_2\, \tau}{\tau+ \sigma_2(x)}
\,,\,
\frac{\nu_3\,\tau }{\tau+ \sigma_3(x)}
\Big)\, .
$$
\end{definition}

In the next discussion this is used with $\nu$ equal to the outward unit normal to 
$\partial\caQ_\delta$.  Next choose a boundary
condition for the stretched equations on $\caQ_\delta$.
 On the flat parts of  $\partial\caQ_\delta$ one has
  $u\in \caE^+(\nu)$.  On the curved parts of the boundary and for 
  $\tau>0$ and real, the stretched problem is symmetric hyperbolic.
  \JR{The normal is $\nu$ and the coefficient of $\partial_j$ is  $A_j\tau/(\tau +\sigma_j)$
  so 
  the normal matrix is $A(\widetilde\nu(\tau,x))$. The maximally
  dissipative condition is 
  $u\in \caE^+(  \widetilde\nu)$.   }
  If $u(\tau)$ is holomorphic and satisfies this condition for $\tau>0$ then 
  by analytic continuation it 
    holds for general $\tau$.   Therefore,  $u\in \caE^+(\widetilde \nu)$
    is a natural maximally dissipative condition for $\tau$ complex.

 The main  result for the stretched  system on 
 $\caQ_\delta$ is the following.  
 
 \begin{theorem}
\label{thm:deltastretched}
For $0<\ell<1$ there exist $C,M_1$ so that for all $\delta\in(0,1)$, 
$M\ge M_1$,  and  holomorphic
$F:
\{ { \rm Re}\, \tau>M\} \to  C^\infty_{\ell\overline\Q}(\caQ_\delta)$, 
there is a unique holomorphic 
$u^\delta:
\{ { \rm Re}\, \tau>M\}
\to H^2(\caQ_\delta)$
 satisfying
\begin{equation}
\label{eq:stretchedpair}
\begin{aligned}
L(\tau,\widetilde\partial_x)  {u^\delta}\ =\ F,
\quad
{\rm on}
\quad \caQ_\delta,
\qquad
{u^\delta}|_{\partial\caQ_\delta} \,\in \, \caE^+\big(A(\widetilde\nu(\tau,x))\big).
\end{aligned} 
\end{equation}
In addition, 
\begin{equation}
\label{eq:elliot}
\begin{aligned}
({\rm Re}\, \tau)\,
\big\| {u^\delta}
&\big\|_{L^2(\caQ_\delta)}
\ +\ 
({\rm Re}\, \tau)^{1/2}\,
\big\| {u^\delta}  
  \big\|_{L^2(\partial\caQ_\delta)}
\cr
\  &+\ 
  \frac{\Re\tau}{|\tau|}\,
  \|\nabla_x u^\delta\|_{L^2(\Q_\delta)}
\ \le\ 
C\,  \big\|
F(\tau)
\big\|_{L^2_{\ell\overline\Q}(\caQ_\delta)}
.
\end{aligned}
\end{equation}
\end{theorem}

{\bf Strategy of proof.}   
 \JR{Theorem \ref{thm:deltastretched} will be  proved
 in 
 Section
 \ref{sec:stretcheddelta}   
 by    
 solving
carefully constructed Helmholtz equations
and boundary conditions on 
$\caQ_\delta$.  
On
$\caQ_\delta$ the solutions are smooth. 
The smoothness is used to 
 prove that the solution of the Helmholtz problem
on $\caQ_\delta$ solves the stretched Pauli system when
proving Theorem \ref{thm:deltastretched}.
Taking the limit $\delta\to 0$  in Section \ref{sec:stretched}    yields 
Theorem \ref{thm:mainstretched}.  }

\subsection{Second boundary condition for the Helmholtz BVP}
\label{sec:secondBC}

\JR{ Theorem
\ref{thm:deltastretched} 
concerns a boundary value problem for the
stretched Pauli system.
One starts from the stretched Pauli system and 
the single boundary condition
$u\in \caE^+(\widetilde\nu)$.  

The stretched Pauli system implies a  stretched 
Helmholtz system of second order.   That $2\times2$ system of second order 
requires 
two  boundary
conditions.

Corollary   \ref{cor:BC2} of this section yields a  crucial
second  boundary
condition.  Example \ref{ex:weakform} shows that
it is a natural boundary
condition for a weak  formulation. 

Section \ref{sec:stretcheddelta} includes a proof of 
the converse implication 
that the 
Helmholtz equation plus the two boundary conditions
imply the  stretched Pauli equations.}

\subsubsection{Neumann identity for the unstretched Pauli system}
\label{sec:unstretchedPauli}

\JR{This section proves that  at the boundary the $\pi^+$ projection of 
the operator $\sum A_j\partial_j$ is a Neumann type boundary operator.
This is used to generate the second boundary condition that is needed
to construct a boundary value problem for the Helmholtz system introduced
in the following 
sections.
}

\begin{definition}   
\label{def:Weingarten}
  For $\ux\in \partial{Q_\delta}$ the {\bf Weingarten map}
  (see for example \cite{hicks})
   is 
 the real selfadjoint map of the 
  tangent space $T_\ux(\partial{Q_\delta})$ to itself that is the 
  differential of the unit exterior normal $\nu$. It  maps
   \mbox{$
 T_\ux(\partial{Q_\delta}) \ni  \bfv
 \to\
    \bfv\cdot\nabla \nu$}.
  Its
  eigenvalues are the {\bf principal curvatures} of $\partial{Q_\delta}$ at $\ux$.
 The {\bf mean  curvature},  
  denoted
$H_{Q_\delta} (\ux)$, is the 
  average of the two principal curvatures.
 
 Extend $\nu$ to a smooth unit vector field defined on a neighborhood
 of $\partial{Q_\delta}$ so as to be constant on  normal lines to the boundary.
 Then $\pi^\pm (\nu(x))$ is well defined
 and smooth  for $x$ in a neighborhood of $\partial{Q_\delta}$. 
\end{definition}

The term $H_{\Q_\delta}(\ux)$ is equal to zero except for a $\delta$ neighborhood of
$\caS$ where it attains values $\sim 1/\delta$.  
The identity of the next proposition is simple in the case of 
flat boundaries.

 \begin{proposition}{\JR{\bf (Neumann identity)}}
 \label{prop:atlast4}  If 
 $u\in H^2({Q_\delta})$ satisfies the boundary  condition
 $
 \pi^-\big( \nu\big) \, u \, =\, 0$
  {\rm on}
 $
 \partial{Q_\delta},
 $
 then,
 \begin{equation}
\label{eq:neumannbc7}
\pi^+\big(\nu\big)
\sum_{j=1}^3 A_j \partial_j u
\ =\ 
\pi^+(\nu) 
\big(\nu\cdot \partial_x +  2H_{Q_\delta} \big)
u,
\quad
{\rm on}
\ \ \partial{Q_\delta}
 .
 \end{equation}

  \end{proposition}

  {\bf Proof of Proposition \ref{prop:atlast4}.}  
      An 
 invariance argument shows that it is sufficient to treat the case
 where $\ux=0$, $\nu(\ux)=(-1,0,0)$ and 
 the $x_j$-axes for $j\ge 2$ are principal curvature directions
 of $\partial{Q_\delta}$. 
  
 Denote by  $\bfe_j$, $j=1,2,3$
  the standard basis for $\RR^3$.
 The principal curvatures corresponding to the tangent directions
 $\bfe_2$ and $\bfe_3$ are denoted 
  $\kappa_2$ and $\kappa_3$.
  The mean curvature 
  is
 $H:=(\kappa_2 +\kappa_3)/2$.
    At $\ux$ the outward unit normal is $-\bfe_1$.
  At $\ux$ the 
 principal curvature  formulas are $\partial_2\nu = -\kappa_2 \bfe_1$ and 
 $\partial_3\nu= -\kappa_3 \bfe_1$.

\vskip.1cm
{\bf First simplifications of the left hand side of \eqref{eq:neumannbc7}.}
The operator on the left is $\pi^+(\nu)(A_1\partial_1 + A_2\partial_2+A_3\partial_3)$.
On the $x_1$ axis, 
$\nu=(-1,0,0)$, so 
 $\pi^+(\nu(x))A_1 = -\pi^+(\nu(x))$.
 On that axis the operator is 
  \begin{equation}
\label{eq:firstterm2}
- \pi^+(\nu) \partial_1 
  +
   \pi^+(\nu)\big(
A_2\partial_2+A_3\partial_3
 \big)
 =
 \pi^+(\nu(x))\nu(x)\cdot\partial_x
+
\pi^+(\nu)\big(
A_2\partial_2+A_3\partial_3
 \big)
.
 \end{equation}

 \vskip.1cm
{\bf Second simplifications.}  Consider the two summands
$\pi^+(\nu)A_j\partial_ju$ with $j\ge 2$.
   On the $x_1$-axis, part {\bf ii} of Proposition \ref{prop:paulispec2} implies that 
   \begin{equation}
   \label{eq:laurence2}
 \pi^+(\nu)\,A_2\,\pi^+(\nu)
  \ =\
 \pi^+ (\nu)
 \,A_3\,\pi^+(\nu)
  \ =\ 0\,.
  \end{equation}  
 Using the boundary condition yields
 \begin{equation}
\begin{aligned}
\label{eq:charles}
 \partial_j
u
 =
 \partial_j
 \Big(
  \pi^+(\nu)  u
  + \pi^-(\nu)  
 u\Big)
  =
 \partial_j
 \big(
  \pi^+(\nu)  u\big)
 \ \ 
 {\rm at}\ \ \ux  .
\end{aligned}
\end{equation}

 For $j\in \{2,3\}$ if $Z$ is  a vector field  on a neighborhood
of $\ux$ that is tangent to the boundary and satisfies  $Z(\ux)=\partial_j$
then
$$
\partial_j u(\ux)
\ =\ 
Z\big(u|_{\partial{Q_\delta}}\big)(\ux)\,.
$$
Since $\pi^+(\nu) u =u$ on the boundary it follows that
$$
\partial_ju(\ux) \ =\ 
Z\big( \pi^+(\nu)    u|_{\partial{Q_\delta}}\big)(\ux)
\ =\
\Big(\partial_j
\big(
\pi^+(\nu)   u\big)\Big)(\ux)\,.
$$
Using \eqref{eq:laurence2} in the last of the following equalities yields
\begin{equation}
\label{eq:steph2}
\begin{aligned}
 \pi^+   (\nu )  
 \,
  A_j 
 \partial_j
u(\ux)
\  &=\ 
   \pi^+(\nu)
   \,
A_j
\Big(
 \partial_j
 \big[
  \pi^+(\nu)
 u\big] \Big)
 \big(\ux\big)\
 \cr
 & =\
    \pi^+(\nu)
   \,
A_j
\Big(
\partial_j \pi^+(\nu )\,u(\ux)
\, +\,
\pi^+ (\nu)\,\partial_ju(\ux)
\Big)
\cr
& =\ 
 \pi^+(\nu)
 \,
A_j\big(
\partial_j \pi^+(\nu) 
\big)
u(\ux)
 \,.
 \end{aligned}
\end{equation}

\vskip.1cm
{\bf The perturbation theory step.}
Use perturbation theory to compute the term  $\partial_j\pi^+$
in the last  expression.
Denote by $Q(\xi)$  the partial inverse of $A(\xi)-|\xi|I$ 
associated to the eigenvalue $+|\xi|$.  It is defined
by 
$$
Q(\xi)\big(
A(\xi)-|\xi|I\big)\,=\,
I-\pi^+(\xi), \qquad
Q(\xi)\,
\pi^+(\xi)
\,=\, 0\,.
$$
Writing 
$$
A(\xi) - |\xi|\,I = \big( |\xi|\pi^+ - |\xi| \pi^-) 
\ -\ 
\big(|\xi|\pi^+ + |\xi|\pi^-)
\ =\ 
-2|\xi|\pi^-
$$
 shows that $Q = (-2|\xi  |)^{-1}\pi^-(\xi)$.

First order perturbation theory (
 Theorem 3.I.2 in \cite{Rauch:LGO},
or 
formulas (II.2.13), (II.2.33)  in 
\cite{kato})
implies that 
\begin{equation}
\label{eq:pertformula}
\frac{\partial}{\partial x_j}
 \Big(   \pi^+\big(A(\nu)\big)\Big) =
-\, 
\pi^+(\nu)\,
\Big(\frac{\partial  A(\nu)}
{\partial x_j} 
\Big)
Q(\nu) 
 -
Q(\nu) 
\Big(
\frac{\partial  A(\nu)}
{\partial x_j} 
\Big)
\pi^+(\nu).
\end{equation}

\vskip.1cm
{\bf  Endgame.}
When 
\eqref{eq:pertformula}  is injected in \eqref{eq:steph2} the
contribution of the  first term vanishes thanks to 
\eqref{eq:laurence2}. 
Turn next to 
$$
\frac{\partial}
{\partial x_j} 
A(\nu(x))
\ =\
A\Big(
\frac{\partial\nu}{\partial  x_j}
\Big)\,.
$$
The principal curvature  formulas imply that at $\ux$,
$$
\frac{\partial\nu}{\partial  x_j}
\ =\ 
\kappa_j(\ux)\,
\bfe_j\,,
\quad
{\rm for}
\quad
j=2,3\,,
\qquad
{\rm so},
\qquad
A\Big(
\frac{\partial\nu}{\partial  x_j}
\Big)\ =\ 
\kappa_j(\ux)\,A_j\,.
$$

Therefore \eqref{eq:steph2} yields
$$
 \pi^+(\nu)\,
   A_j
 \partial_j
u(\ux)
\ =\
\kappa_j(\ux)\,
 \pi^+(\nu)\,
A_j
\,
\pi^-(\nu)
A_j
\,
\pi^+(\nu)\,.
$$
Compute using
\eqref{eq:laurence2} and omitting the argument $\nu(\ux)$ for ease of reading yields
\begin{equation*}
\pi^+A_j\pi^-A_j\pi^+
=
\pi^+A_j\big(\pi^- + \pi^+\big)A_j\pi^+
=
\pi^+A_j\,A_j\pi^+
=\pi^+\pi^+
=
\pi^+\,.
\end{equation*}

Therefore
 \begin{equation}
 \label{eq:j23}
 \pi^+(\nu(\ux))\,
   A_j
 \partial_j
u(\ux)
 =
\kappa_j(\ux)\,
 \pi^+(\nu(\ux)) \,u,
 \quad
 {\rm for}
 \quad
 j=2,3.
 \end{equation}
  
   The sum of the terms \eqref{eq:j23}   is equal to
 $(\kappa_2 +\kappa_3)\, \pi^+u=2H_{Q_\delta}\,\pi^+u$.  This
 yields 
\begin{equation*}
\label{eq:neumannbc2}
\pi^+(\nu(\ux))\Big(\nu\cdot \nabla_x +  2H_{Q_\delta}(\ux)\Big)
 u
  \,.
 \end{equation*}
This completes the proof of \eqref{eq:neumannbc7}.
  \hfill
 $\Box$
 \vskip.2cm

 \subsubsection{Transverse identity  for stretched Pauli for  $\tau\in]m,\infty[$}
 \label{sec:stretchedPauli}

 \JR{ Recall that  for real $\tau$, $\uOmega$ is the image of ${Q_\delta}$
  by the stretching transformations in Definition \ref{def:Xj}. }
 To find  the conormals to  $\uOmega$, compute
  $$
  \sum_j \nu_j dx_j  =
  \sum_j \nu_j \sum_k \frac{\partial x_j}{\partial X_k } dX_k
   =
  \sum_j
  \nu_j
  \frac{\partial x_j}{\partial X_j } dX_j
   =
\sum_j
  \frac{  \nu_j\, \tau}{\tau + \sigma_j}\,
  \, dX_j\, .
  $$
  $  \sum_j \nu_j dx_j$
   annihilates the tangent space  to $\partial{Q_\delta}$ at $x$.
 The map $x\to X$ takes the tangent  space to ${Q_\delta}$ to the 
 tangent space to $\uOmega$,  
 Therefore,
 $\sum_j
   \nu_j\tau/( \tau + \sigma_j)
  \, dX_j$ 
  annihilates the tangent space to  $\uOmega$
  at $X(x)$.  
 It is therefore  a conormal to $\uOmega$.    The unit conormal $\nu_{\uOmega}(X)$  is 
  $$
 \nu_\uOmega(X)
 \ =\
  \bigg(\sum_j
    \frac{ \nu_j^2(x(X))\, \tau^2  }
   {(\tau+\sigma_j(x(X)))^2}\bigg)^{-1/2}
  \bigg(
  \sum_{j=1}^3
  \frac{ \nu_j(x(X) ) \,\tau   }
   {\tau+\sigma_j(x(X))}\,dX_j
   \bigg)
   \,.
   $$

  \begin{definition}
  \label{def:V}
  For ${\rm Re}\,\tau>0$ 
   and $x$ on a neighborhood of $\partial{Q_\delta}$,
  define the
  first order  differential operator
 $V$
 by 
 \begin{equation}
 \label{eq:defV}
 \begin{aligned}
 V(\tau, x,\partial)\ &:=\ 
  \Big(\sum_j
    \frac{\nu_j^2\,\tau^2}
   {(\tau+\sigma_j)^2}\Big)^{-1/2}\
\sum_j
  \frac{ \nu_j\,\tau^2}
   {(\tau+\sigma_j)^2}\,
   \frac{\partial}{\partial x_j}\,.
   \end{aligned}
\end{equation}
\end{definition}

\begin{remark}
{\bf i.}   For $\tau\in ]0,\infty[$,
$V$ is a unit vector  field transverse to $\partial{Q_\delta}$
since its scalar product with the unit outward normal $\nu\cdot\partial$ is strictly positive.

{\bf ii.}  For  $\tau$ not real, the coefficients of 
 $V$  are not real, so $V$ is not
 a vector field.  
 
{\bf iii.}   There is an $R>0$ independent of 
 $\delta$ so that for $|\tau|>R$, $\partial{Q_\delta}$ is non characteristic for $V$.
 {\rm
 Indeed,  $V-\nu\cdot \partial$ has coefficients $O(1/\tau)$ and 
 the boundary is noncharacteristic for 
 $\nu\cdot\partial$.  
 }
 \end{remark}

  \begin{corollary} \JR{
  \bf  (Transverse identity 1)
  }
  
 \label{cor:neumannrealtau}  There is an $m>0$ so that if 
 $\tau\in ]m,\infty[$ and
  $u\in H^2({Q_\delta})$ satisfies the boundary
  condition
 \begin{equation}
 \label{eq:paulitilde}
 u\ \in\
 \caE^+\big(\widetilde \nu(\tau,x ) \big)
   \quad
   {\rm on}\quad
   \partial{Q_\delta},
\end{equation}
then with   $V(\tau,x,\partial)$ from \eqref{eq:defV},  
\begin{equation}
\label{eq:neumannbc4}
 \pi^+ 
 (\widetilde \nu  )\,
\sum_{j=1}^3 A_j \widetilde\partial_j u 
\ =\ 
  \pi^+
 (\widetilde \nu )\,
  \Big(V(\tau,\ux,\partial)+ 2H_\uOmega(X(\tau,\ux))
  \Big)u
  \ \ 
  {\rm on}
  \ \ \partial{Q_\delta}\,.
  \end{equation}
 
  \end{corollary}
 
 \begin{remark}
  The normal matrix of the stretched system is equal to 
  $A(\widetilde \nu)$.
  For positive $\tau$, the boundary condition in
  \eqref{eq:paulitilde} is the natural
 maximally absorbing one.  
  \end{remark}

  {\bf Proof of Corollary \ref{cor:neumannrealtau}.}        Define $v:\uOmega\to \CC^2$ by 
  $v(X):= u(x(X))$.   
 Since $u$ satisfies the stretched
  Pauli system on a neighborhood of $\partial{Q_\delta}$, 
  \eqref{eq:cov} implies that 
    $v$ satisfies the unstretched Pauli system
  on a neighborhood of $\partial\uOmega$.
  
   The unstretched differential equation satisfied by $v$ has 
   principle symbol  
   $
   \sum_j A_j\, {\partial}/{\partial X_j}\,.
   $
  The symbol at any outward  conormal vector to $\uOmega$ is equal to a
  positive multiple of 
  $\sum_j A_j \nu_j \tau/(\tau + \sigma_j)$.  This 
  sum is equal to the symbol of the 
  stretched operator on ${Q_\delta}$ at the conormal $\nu$
  to ${Q_\delta}$. Thus
  the positive eigenspace of the unstretched symbol at $\nu_\uOmega(X)$
  is equal to the positive eigenspace of the stretched operator at 
  $\nu_{Q_\delta}(x)$.  
  
 The boundary condition satisfied by $u$ asserts that 
 $$
 u\ \in\
 \caE^+\big(
  A(\widetilde \nu) \big)
    \ =\
    \caE^+\big(
  A\big(\nu_\uOmega)\big)
   \,.
   $$
Therefore
$v$ satisfies the boundary condition
$
v |_{\partial\uOmega}\in
  \caE^+\big(
  A(\nu_\uOmega)
\big)$.
The function $v$ on $\uOmega$ therefore satisfies the hypotheses
of 
  Proposition \ref{prop:atlast4} on $\uOmega$.
  That Proposition  implies that
for $X\in \partial\uOmega$,
  \begin{equation*}
    \pi^+
  \big(  
 \widetilde \nu(x(X))\big)  
\sum_{j=1}^3 A_j \widetilde\partial_j u 
\ =\
 \pi^+
  \big(
 \widetilde \nu(x(X)  ) \big)
 \Big(
   \nu_\uOmega\cdot \partial_X 
   + 2H_{\uOmega}(X)
   \Big)
  v\,.
   \end{equation*}
 Equation \eqref{eq:cov}   shows that 
    \begin{equation*}
  \begin{aligned}
  \nu_\uOmega\cdot \partial_X
  \ =\ 
   \Big(\sum_j
    \frac{ \nu_j^2}
   {(\tau+\sigma_j)^2}\Big)^{-1/2}
\sum_j
  \frac{ \nu_j}
   {\tau+\sigma_j}
   \
   \frac {\tau}{\tau+\sigma_j}
   \,
   \frac{\partial}{\partial x_j}
   \ =\ V\,.
   \end{aligned}
     \end{equation*}
  Inserting in the preceding equation yields
   \eqref{eq:neumannbc4}.
   \hfill
   $\Box$
   \vskip.2cm

    \subsubsection{Transverse identity for stretched Pauli for $\tau\notin\RR$ }
 \label{sec:stretchedPauli2}
 
 \JR{The next proposition shows that several quantities depend
 holomorphically on $\tau$.  
 Part {\bf iv} of the  next proposition is
the key identity  for complex $\tau$.   It follows from the real
identity by analytic continuation.
}
 
 \begin{proposition}  \label{prop:analyticBC}
   {\bf i.}  There is  an $R_1>1$ so that for $|\tau|>R_1$ the spectrum of 
 $A(\widetilde \nu(\tau,x))$
 consists of one simple eigenvalue  in $|z-1|<1$ and a second
in $|z-(-1)|<1$.      
Then 
 the map $\tau \mapsto \pi^\pm\big(A\big(
\widetilde \nu(\tau,x)\big)$
 is analytic in 
$|\tau|>R_1$.   

 \vskip.1cm
 
 {\bf ii.}        There is an $R_2\ge R_1$ so that the function  $\tau\mapsto \nu_\uOmega(X(\tau,x))$ from
 $]m,\infty[$ to $C^\infty(\partial{Q_\delta})$ has a holomorphic extension to
 $\{|\tau|>R_2\}$.

 \vskip.1cm
 
 {\bf iii.}   There is an $R_3\ge R_2$ so that the function  $\tau\mapsto H_{\uOmega}(X(\tau,x))$ from
 $]m,\infty[$ to $C^\infty(\partial{Q_\delta})$ has a holomorphic extension to
 $\{|\tau|>R_3\}$.

\vskip.1cm
{\bf iv.}   \JR{
 \bf
 (Transverse identity 2)
 }
If $|\tau|>R_3$ and
 $\tau \mapsto u(\tau)\in H^2({Q_\delta})$ satisfies $u\in \caE^+(\widetilde\nu)$ on $\partial{Q_\delta}$
 and is holomorphic on a connected open subset 
 $ {\cal O}   \subset \{\tau\in \CC:  |\tau|>R_3\}$ that meets the real axis in 
 a nonempty open set, 
   then
 \eqref{eq:neumannbc4} holds on $\caO$.

\end{proposition}

{\bf Proof.}  {\bf i.} For $|\tau|$ large one has uniformly for $x\in \partial{Q_\delta}$,
$$
\widetilde \nu \ =\ 
\bigg(
  \frac{ \nu_1\,\tau}
   {\tau+\sigma_1}
\,,\,
\frac{ \nu_2\,\tau}
   {\tau+\sigma_2}
  \,,\,
\frac{ \nu_3\,\tau}
   {\tau+\sigma_3}
   \bigg)
   \ =\ \nu \,+\, O(|\tau|^{-1})\,,
   $$
   The assertions
   in  {\bf i} follows from Part {\bf ii} of Proposition \ref{prop:paulispec}.
   
   \vskip.1cm

{\bf ii.} It suffices to construct the analytic continuation for points in  a neighborhood of each
$\uX\in \partial\uOmega$.  
Suppose that $\uX=X(\tau,\ux)$ with $\tau>0$ and 
$\ux\in \partial{Q_\delta}$  
and $X$ the stretching transformation defined 
by 
\eqref{eq:defXj}.  
The map  
$\tau\mapsto X(\tau,\cdot)$ is holomorphic on $\tau\ne 0$ with values
in $C^\infty(\partial{Q_\delta};\CC)$.
In addition, 
$\partial X/\partial x = I + O(1/\tau)$, so
$\partial X/\partial x$  is invertible for 
$|\tau|>R$.   

Suppose that 
$x(\alpha_1,\alpha_2)$ is a parametrization of a neighborhood of $\ux$ in $\partial{Q_\delta}$.
Then for $\tau>0$, $X(\tau, x(\alpha_1,\alpha_2))$ is a parametrization of a neighborhood of 
$\uX$ in $\partial\uOmega$.   
For those $\tau$ the tangent space to $\partial\uOmega$ is spanned by
the independent vectors $\partial X(\tau,x(\alpha))/\partial \alpha_i$, $1=1,2$. 
Thanks to the invertibility of $\partial X/\partial x$, the formula
\begin{equation}
\label{eq:TX}
{\rm Span}\  \Big\{
\frac{\partial X(\tau,x(\alpha)) }{ \partial \alpha_1}
\ ,\
\frac{\partial X(\tau,x(\alpha)) } { \partial \alpha_2}
\Big\}
\ =\
{\rm Span}\  \Big\{
\frac{\partial X}{\partial x}
\frac{\partial x}{\partial \alpha_1}
\ ,\
\frac{\partial X}{\partial x}
\frac{\partial x}{\partial \alpha_2}
\Big\}
\end{equation}
shows that the tangent space has a holomorphic continuation to $| \tau|>R$.

For real $\tau$ a normal vector to $\uOmega$ at $X(\tau,x(\alpha))$ is given by
$$
\frac{\partial X}{\partial x}
\frac{\partial x}{\partial \alpha_1}
\ \wedge\
\frac{\partial X}{\partial x}
\frac{\partial x}{\partial \alpha_2}.
$$
It is nonvanishing because $\partial X /\partial x$ is invertible and
the vectors $\partial x/\partial \alpha_j$ are independent.
The unit normal vector is given by
$$
\nu(X(\tau,x(\alpha))\ =\ 
\frac{
\frac{\partial X}{\partial x}
\frac{\partial x}{\partial \alpha_1}
\ \wedge\
\frac{\partial X}{\partial x}
\frac{\partial x}{\partial \alpha_2}
}
{
\big[
\sum_i
\big(\big(
\frac{\partial X}{\partial x}
\frac{\partial x}{\partial \alpha_1}
\ \wedge\
\frac{\partial X}{\partial x}
\frac{\partial x}{\partial \alpha_2}
\big)_i\big)^2
\big]^{1/2}
}
$$
Since $\partial X/\partial x = I + O(1/\tau)$ it follows that one can
choose $R>0$ so that
  $$
  \sum_i
\bigg(\Big(
\frac{\partial X}{\partial x}
\frac{\partial x}{\partial \alpha_1}
\ \wedge\
\frac{\partial X}{\partial x}
\frac{\partial x}{\partial \alpha_2}
\Big)_i\bigg)^2
$$
has strictly positive real part for $|\tau|>R$.   With that choice the expression
for $\nu(X(\tau,x(\alpha))$
yields an analytic continuation of the unit normal vector to $|\tau|>R$.  For
$\tau\notin \RR$,
$\nu(X(\tau,x(\alpha))$ need  not be real and need  not be of unit length.

{\bf iii.}  For $\tau$ real the Weingarten map is the map from $T_\uX(\partial\uOmega)$
to itself that maps the two basis vectors as follows,
\begin{equation}
\label{eq:weingarten}
\frac{\partial X(\tau, x(\alpha))}{\partial \alpha_j}
\quad
\to
\quad 
\frac{\partial \nu_\uOmega(X(\tau, x(\alpha))}
{\partial \alpha_j},
\qquad
j=1,2.
\end{equation}
The holomorphic extension of $\nu$ implies that 
 the Weingarten map extends holomorphically to a family of
linear map of the holomorphic family of  two dimensional spaces
\eqref{eq:TX}
to itself.     

For $\tau$  real the mean curvature  $H_\uOmega$ is equal to one half of the trace
of  the Weingarten map.   The preceding paragraph shows that 
this trace has a holomorphic continuation proving {\bf iii.}

{\bf iv.}   
The difference of the two sides  of  \eqref{eq:neumannbc4} is holomorphic on 
a connected set.
Corollary \ref{cor:neumannrealtau} implies that it vanishes for $\tau$
on the open intersection with the real axis.   By analytic continuation 
it vanishes identically.
\hfill
\qed
\vskip.2cm

\begin{corollary} 
\label{cor:BC2}
If $u\in H^{2}({Q_\delta})$ satisfies $L(\tau, \widetilde \partial)u=0$ 
on $\partial{Q_\delta}$ and 
$u\in \caE^+(\widetilde\nu)$ on $\partial{Q_\delta}$, 
then
\begin{equation}
\label{eq:BC2}
 \pi^+
 (\widetilde \nu )
  \Big(V(\tau,x,\partial)+ \tau+2H_\uOmega(X(\tau,  \ux)) \Big) u=0
 \quad
  {\rm on}
  \quad
  \partial{Q_\delta}\,.
\end{equation}
\end{corollary}

{\bf Proof.}  
Equation
\eqref{eq:neumannbc4} 
implies that 
\begin{equation*}
\pi^+(\widetilde \nu) L(\tau, \widetilde \partial) u\ =\ 
 \pi^+
 (\widetilde \nu )\,
  \Big(V(\tau,x,\partial)+ \tau+2H_\uOmega(X(\tau,  \ux))
  \Big)u  
  \quad
  {\rm on}
  \quad
  \partial{Q_\delta}\,.
\end{equation*}
Since 
$L(\tau, \widetilde \partial)u=0$ on   $\partial{Q_\delta}$,
$u$  satisfies
\eqref{eq:BC2}.
\hfill
\qed
\vskip.2cm

\section{Analysis of the Helmholtz BVP}
\label{sec:PauliHelmholtz}

This section derives Helmholtz equations in two steps.
The first repeats the usual derivation of Helmholtz from 
Pauli adapted to the stretched operators.
It yields an
 operator that is not in divergence 
 form.

For $i\ne j$ the anticommutation formulas  \eqref{eq:anticommutation}   
imply that 
$$
A_i\widetilde \partial_i
\,
A_j\widetilde \partial_j
\ +\ 
A_j\widetilde \partial_j
\,
A_i\widetilde \partial_i
\ =\ 0,
\qquad
{\rm for}
\qquad
i\ne j\,.
$$ 
Indeed,  when the derivatives 
fall on variable coefficients they yield zero.
Define
$$
\widetilde\partial_j^2 
\ :=\ 
\Big(
\frac{\tau}{\tau + \sigma_j(x_j)}
\,
\partial_j
\Big)
\Big(
\frac{\tau}{\tau + \sigma_j(x_j)}
\,
\partial_j
\Big)
$$
where  the order of the operators 
inside the parentheses
is important.
The following stretched versions of 
\eqref{eq:anticomm1} hold,
\begin{equation}
\label{eq:related}
\begin{aligned}
&\Big(
\sum_j A_j\widetilde \partial_j\Big)^2
\  =\ -\,
\sum_j \widetilde\partial_j^2 \,,
\cr
&\Big(\sum A_j\widetilde\partial_j-\tau\Big)
\Big(\sum A_j\widetilde\partial_j+\tau\Big)
\  =\
 \sum_j \widetilde\partial_j^2
 \ -\
 \tau^2\,.
 \end{aligned}
\end{equation}

The second equation in \eqref{eq:related}  shows that where
a function 
$u$ satisfies  $L(\tau, \widetilde \partial) u=0$,
it satisfies 
$\big( \sum_j \widetilde\partial_j^2 -\tau^2\big)u=0$

\subsection{The  Helmholtz identity}
\label{sec:helmop}

\JR{The  next lemma shows that multiplying 
by a suitable weight yields an operator in divergence
form.   This is used in the derivation of {\it a priori} estimates.}

\begin{definition}  
\label{def:p}
\JR{
In the next formula,  when the subscript does not belong to $\{1,2,3\}$,
it is replaced by the unique element of that set that is congruent
modulo $3$.  
Define the scalar divergence form operator}
 \begin{equation}
\label{eq:tildewave2p}
p \big(
\tau, x, \partial \big)u
\ :=\
\sum_{j=1}^3\ \partial_j\ 
\frac{(\tau + \sigma_{j+1}(x_{j+1})) (\tau+\sigma_{j+2}(x_{j+2} ))}{\tau ( \tau + \sigma_j (x_j)) }
\ 
\partial_j u \,,
\end{equation}
and
\begin{equation}
\label{eq:defPi}
\Pi(\tau,x)\ :=\
\prod_{i=1}^3\ 
\frac{\tau + \sigma_i(x_i)}{\tau}\,.
\end{equation}
\end{definition}

\begin{lemma}
\JR{
\bf 
(Stretched Helmholtz identity)
}
\label{lem:helmder}
As operators on $H^2_{loc}(\caQ)$,
\begin{equation}
\label{eq:prodpauli}
\Pi(\tau,x)
\Big(\sum_j A_j\widetilde\partial_j   -\tau\Big)
\Big(\sum_j A_j\widetilde\partial_j+\tau\Big)
\, =\,
\Big(
p(\tau,x,\partial) 
\, -\, 
\tau^2\,\Pi(\tau,x)
\Big)I
.
\end{equation}

\end{lemma}

{\bf Proof.}
Expanding the product on the left
using the anticommutation relations \eqref{eq:anticommutation}  yields
$$
\sum_j\
\bigg(
\prod_{i=1}^3
\frac{\tau + \sigma_i(x_i)}{\tau}
\bigg)
\
\Big(
\frac{\tau}{\tau + \sigma_j(x_j)}
\,
\partial_j
\Big)
\Big(
\frac{\tau}{\tau + \sigma_j(x_j)}
\,
\partial_j
\Big)
\ -\ 
\tau^2\,
\prod_{i=1}^3
\frac{\tau + \sigma_i(x_i)}{\tau}
\,.
$$
The factor before the first derivative on the left is equal to 
$$
\bigg(
\prod_{i=1}^3
\frac{\tau + \sigma_i(x_i)}{\tau}
\bigg)
\
\frac{\tau}{\tau + \sigma_j(x_j)}
\ =\
\frac{
\big(
\tau + \sigma_{j+1}(x_{j+1})\big)
\big(
\tau + \sigma_{j+2}(x_{j+2})\big)
}
{\tau^2
}
\,.
$$
This function does not depend on $x_j$ so commutes with 
$\partial_j$.
\begin{equation*}
\begin{aligned}
\bigg(\prod_{i=1}^3
\frac{\tau + \sigma_i(x_i)}{\tau}\bigg)
\ &
\Big(
\frac{\tau}{\tau + \sigma_j(x_j)}
\,
\partial_j
\Big)
\Big(
\frac{\tau}{\tau + \sigma_j(x_j)}
\,
\partial_j
\Big)
\cr
&=\
\partial_j
\Big(
\frac{
\big(
\tau + \sigma_{j+1}(x_{j+1})\big)
\big(
\tau + \sigma_{j+2}(x_{j+2})\big)
}
{\tau^2
}
\
\frac{\tau}
{\tau + \sigma_j(x_j)}
\Big)
\partial_j
\cr
&=\
\partial_j
\Big(
\frac{
\big(
\tau + \sigma_{j+1}(x_{j+1})\big)
\big(
\tau + \sigma_{j+2}(x_{j+2})\big)
}
{\tau\,
\big(\tau + \sigma_j(x_j)\big)}
\Big)
\partial_j\,.
\end{aligned}
\end{equation*}
This completes the proof.
\hfill
$\Box$
\vskip.2cm

\begin{remark}
{\bf i.}  The factors in the product on the left of \eqref{eq:prodpauli} are
\begin{equation}
\sum_j A_j\widetilde\partial_j+\tau
\ =\ 
L\big(\tau, \widetilde\partial\big)\,,
\quad
{\rm and},
\quad
\sum_j A_j\widetilde\partial_j   -\tau
\ =\ 
L\big(-\tau, \widetilde\partial\big)\,.
\end{equation}

{\bf ii.}  Since
\begin{equation}
\label{eq:nearlyH2}
\Big|
\Pi(\tau,x) -1
\Big|
\ =\ 
\Big|\prod_{i=1}^3
\frac{
\tau + \sigma_i(x_i)}{\tau}\ -1\Big|
\ \lesssim\ \frac1{|\tau|}
\end{equation}
 the coefficients of the 
 operator on the right of 
  \eqref{eq:prodpauli}  differ
  from those of 
   the classical Helmholtz operator
   $\Delta-\tau^2$ by $O(|\tau|^{-1})$.
\end{remark}

\begin{definition}
$\bullet$  For vectors $\alpha,\beta$ in $\CC^k$ define
$\alpha\cdot \beta :=\sum_j \alpha_j\,\beta_j$.

$\bullet$ Define the continuous bilinear form 
$a: H^1(\caQ;\CC^2)\times H^1(\caQ;\CC^2)
\to \CC$  associated to $\,-\,p$ from Definition \ref{def:p}
by
\begin{equation}
\label{eq:pform}
a(u,v)
  =
 \int_\caQ
\sum_{j=1}^3\  
\frac{(\tau + \sigma_{j+1}(x_{j+1})) (\tau+\sigma_{j+2}(x_{j+2} ))}
{\tau ( \tau + \sigma_j (x_j)) }
\,
\partial_ju
\,\cdot \,
{\partial_j v}
\ dx\,.
\end{equation}

$\bullet$  The formula with integration over
${\caQ_\delta}$ defines  a continuous 
form from $H^1_{loc}({\caQ_\delta})\times H^1_{compact}({\caQ_\delta})
\to \CC$.
\end{definition}

\begin{remark} 
\label{rem:forma}
 {\bf i.}
If $u\in H^1_{loc}({\caQ_\delta})$ 
and $f\in H^{-1}_{loc}({\caQ_\delta})$ then $u$ satisfies $pu=f$   on
${\caQ_\delta}$ if and only if 
$$
\forall \phi\in C^\infty_0({\caQ_\delta}),
\qquad
a(u,\phi) \ =\ -\,
\int_{\caQ_\delta}
 f\cdot \phi
 \ dx\,.
$$ 

\noindent
{\bf ii.}  Multiplying numerator and denominator of the coefficient
of $\partial_j$ 
in \eqref{eq:pform}
by  $\tau + \sigma_j$
shows that 
\begin{equation}
\label{eq:nicea}
\begin{aligned}
a(u,v)
 \ =\ 
 \int_{\caQ_\delta}
\Pi(\tau,x)\ 
   \sum_{j=1}^3
   \frac{\tau^2}{   (\tau +\sigma_j)^2}\
   \partial_ju
\,\cdot\,
{\partial_j v}
\ dx
\,.
\end{aligned}
\end{equation}
\noindent
{\bf iii.}   If $u\in H^2({\caQ_\delta})$, an integration by parts yields
\begin{equation}
\label{eq:helmnormal}
a( u,v)
  =
  -\, \int_{\caQ_\delta}  p\,{u} \,\cdot \, 
  {{v}} \ dx
    +
   \int_{\partial{\caQ_\delta}}
   \Pi(\tau,x)
       \sum_{j=1}^3
   \frac{\nu_j \,\tau^2}{   (\tau +\sigma_j)^2}
   \,\partial_j {u}
   \,\cdot\,
    {v}\, d\Sigma.
\end{equation}

\end{remark}

To solve the stretched equation, start
by  using 
\eqref{eq:prodpauli} to show that any solution must
satisfy the  divergence form   {\bf  Helmholtz equation} 
\begin{equation}
\label{eq:helm}
\begin{aligned}
\Big(
p(\tau,x,\partial) 
-
\tau^2\, \Pi(\tau,x)
\Big)u
\, =\,
\Pi(\tau,x)\
\Big(\sum A_j\widetilde\partial_j   -\tau\Big)
F
.
\end{aligned}
\end{equation}

 \begin{remark}  
 \label{rem:timeharmonic}
 There is an extensive literature on using the PML technology for the 
 solution of time harmonic scattering problems for the  wave
 equation beginning  
 with  Collino-Monk and  Lassas-Somersalo
 \cite{Collino:1998:PML,
 lassas1998existence, 
 lassas2001analysis,
  BRAMBLE2013209,zbMATH06605778}.
 All  depend
 on  choosing $\sigma_j$ constant outside a compact set and then
 relying on an explicit Green's function for the Helmholtz operator
 $\tau^2 - p$ with $\tau=i\omega$ and $x$ outside that compact
 set.   Rellich's  Uniqueness Theorem and the exponential decay of the
  Green's function
 drives the analysis.   
The operator $p$ and the form $a(\cdot,\cdot)$ appear
in those articles.    Variable $\sigma_j$,
corners, and 
absorbing boundary  conditions at trihedral corners have no analogue
in their work.   This time harmonic work is related to the method
of complex scaling in Scattering Theory
introduced by Balslev-Coombes \cite{balslev1971spectral} and raised to high art
by  Sj\"ostrand and a brilliant school (see  \cite{dyatlov2019mathematical}).
\end{remark}
   
  \subsection{Lopatinski for the   Helmholtz BVP}
  \label{sec:HBVP}
  
  \JR{For the equation \eqref{eq:helm} construct a boundary value problem.
  The solutions come from the stretched Pauli system  with the boundary
  condition 
  $u\in \caE^+(\widetilde\nu)$,
  equivalently   
  $\pi^-(\widetilde \nu(\tau,x) ) u =0$.
  Corollary \ref{cor:BC2} 
  provides a  second boundary condition.
  A $2\times 2$ system of second order elliptic
  equations requires exactly two conditions.
  }
 The present section is devoted to studying the
resulting   Helmholtz boundary value problem,
\begin{equation}
\label{eq:PHBVP}
\begin{aligned}
\Big(
\tau^2\,\Pi(\tau,x)\ -\ 
p(\tau,x,\partial) 
\Big)  u \ &=\ f\,
\quad 
{\rm on}
\quad
\caQ_\delta\,,
\cr 
\pi^-(\widetilde \nu(\tau,x) ) u \ &=\ g_1
\quad
{\rm on}
\quad
\partial\caQ_\delta\,,
\cr
\pi^+(\widetilde \nu(\tau,x))
\Big(
V(\tau,x,\partial)+ \tau+2H_\uOmega(X(\tau,\ux))
\Big)u
\ &= \ g_2
\quad
{\rm on}
\quad
\partial\caQ_\delta\,.
\end{aligned}
\end{equation}
Here $g_1$ and $g_2$ are functions on 
$\partial{Q_\delta}$ that take values in 
$\caE^-(A(\widetilde \nu(\tau,x))$ and $\caE^+(A(\widetilde \nu(\tau,x))$
respectively.

\begin{definition}
 For $S\in {\rm Hom}\, \CC^k$ denote by $S^\dagger$ the transposed matrix
  so that 
 $Su\cdot v = u\cdot S^\dagger v$ for all vectors $u,v\in \CC^k$.
 \end{definition}

For $|{\rm Im}\,\xi|<|{\rm Re}\, \xi|$,
$A(\xi)$ has two eigenvalues $\lambda^\pm(\xi)$ and 
spectral representation
 $$
 A(\xi) = \lambda^+\pi^+(\xi)
 +\lambda^-\pi^-(\xi),
 \quad
 {\rm so},
 \quad
 A(\xi)^\dagger= \lambda^+\pi^+(\xi)^\dagger
 +\lambda^-\pi^-(\xi)^\dagger.
 $$
 Therefore $\lambda^\pm$ are eigenvalues of 
 $A(\xi)^\dagger$ and 
 $\pi^\pm(\xi)^\dagger$ are the corresponding spectral projections.

 \begin{definition}
Define the {\bf transposed boundary value problem} as,
\begin{equation}
\label{eq:ADJBVP}
\begin{aligned}
\Big(
\tau^2\,\Pi(\tau,x)\ -\ 
p(\tau,x,\partial) 
\Big)  u \, &=\, f
\quad 
{\rm on}
\quad
\caQ_\delta,
\cr 
\pi^-(\widetilde \nu(\tau,x) )^\dagger u \, &=\, g_1
\quad
{\rm on}
\quad
\partial\caQ_\delta,
\cr
\pi^+(\widetilde \nu(\tau,x))^\dagger
\Big(
V(\tau,x,\partial)+ \tau+2H_\uOmega(X(\tau,\ux))
\Big)u
\, &= \, g_2
\quad
{\rm on}
\quad
\partial\caQ_\delta.
\end{aligned}
\end{equation}
The functions $g_1$ and $g_2$ on $\partial{Q_\delta}$ take values in 
$\caE^-(A(\widetilde\nu(\tau,x))^\dagger)     $
and 
$\caE^+(A(\widetilde\nu(\tau,x))^\dagger)     $ 
respectively.
\end{definition}

In Section
\ref{sec:deltafixed} it is proved that the
 annihilator of the range of the direct problem
is equal to 
the nullspace of the transposed  problem.

\begin{lemma}
There is an $R>0$ independent of $\delta$ to that for $|\tau|>R$
the boundary  
value problems \eqref{eq:PHBVP} and \eqref{eq:ADJBVP} satisfy Lopatinski's condition
for all $x\in \partial\caQ_\delta$.
\end{lemma}

{\bf Proof.}  Analyse only  \eqref{eq:PHBVP}.  The proof for the other
is nearly identical.    

The Lopatinski  condition depends only on the 
highest order terms of the equation and the boundary condition.
Since the highest order terms of $p$ converge to those of $\Delta$
it suffices to verify Lopatinski's condition for $\Delta$ with the 
boundary conditions those of the Helmholtz problem.   This leads to the 
constant coefficient half space problems
\begin{equation}
\label{eq:reducedBVP}
\begin{aligned}
\Delta\, u \ &=\ f\,
\quad 
{\rm on}
\quad
\nu\cdot x<0\,,
\cr 
\pi^-( \nu) u \ &=\ g_1
\quad
{\rm on}
\quad
\nu\cdot x=0\,,
\cr
\pi^+(\nu)
\big(
\nu\cdot \partial_x u\big)
\ &= \ g_2
\quad
{\rm on}
\quad
\nu\cdot x=0\,.
\end{aligned}
\end{equation}
Choosing an orthonormal basis of $\RR^2$ whose first element is a basis
for $\caE^+(\nu)$ and whose second is a basis for $\caE^{-}(\nu)$
reduces to 
\begin{equation}
\label{eq:reducedBVP2}
\begin{aligned}
\Delta\, u \ &=\ f\,
\quad 
{\rm on}
\quad
\nu\cdot x<0\,,
\cr 
u_1 \ &=\ g_1
\quad
{\rm on}
\quad
\nu\cdot x=0\,,
\cr
(\nu\cdot\partial_x)u_2
\ &= \ g_2
\quad
{\rm on}
\quad
\nu\cdot x=0\,.
\end{aligned}
\end{equation}
This is the Dirichlet problem for $u_1$ and the Neumann problem
for $u_2$.  

\JR{The Lopatinski condition concerns the homogeneous problem
with $f=g_1=g_2=0$.   It requires   that for any 
$$
0\ne \xi^\prime\in \RR^d\quad{\rm  with}\quad \xi^\prime\perp\nu\,,
$$
 the only solution
$$
u\ =\ e^{i\xi^\prime\cdot x}\ w(\nu\cdot x),\qquad
w(s)\ \to \ 0 \quad {\rm when}\quad
s\to -\infty\,,
$$
 is 
$w=0$.
For both the Dirichlet and Neumann problems the verification
of the  Lopatinski condition is classical.
}
\hfill
\qed
\vskip.2cm

 \begin{theorem} 
 \label{thm:taufixed} There is an $M>0$  so that if
 ${\rm Re}\,\tau >M$ and $0<\delta<1$, then
 the continuous linear map
 $$
H^2(\caQ_\delta)\ni u\mapsto (f,g_1,g_2)\in
L^2(\caQ_\delta)\times  H^{3/2}(\partial\caQ_\delta;\caE^-(\widetilde\nu))\times
 H^{1/2} (\partial\caQ_\delta;\caE^+(\widetilde\nu))
 $$
 defined by  \eqref{eq:PHBVP}
  is one to one and onto.
  \end{theorem}

  {\bf Strategy of the proof of Theorem \ref{thm:taufixed}.}   The theory of 
  elliptic boundary value  problems  satisfying Lopatinski's condition  implies the following
  facts, see \cite{agmon1959estimates,agmon1964estimates}.

$\bullet$   The kernel of the map is a finite dimensional subset of $C^\infty(\overline{\caQ_\delta})$.

$\bullet$ 
 The range is closed with finite codimension.

 $\bullet$   The annihilator of the range is 
  a subspace of 
 $C^\infty(\overline{\caQ_\delta})
 \times
 C^\infty({\partial\caQ_\delta})
 \times
 C^\infty({\partial\caQ_\delta})
 $.  
 
 {\it To prove the 
 theorem it suffices to prove that  the kernel and the annihilator of the range
 are both trivial.}
 
\JR{
The main step is the proof of a uniform {\it a priori} estimate.  That estimate
is stated in Theorem \ref{thm:perturbed} at the start of the next subsection.  
It's proof is completed at the end of Section \ref{sec:UP}.  With that estimate
in hand, the proof of Theorem \ref{thm:taufixed}  is completed in 
Section \ref{sec:deltafixed}.
}

\subsection{Main {\it a priori} estimate, Theorem \ref{thm:perturbed}}

\begin{theorem} {\bf (Helmholtz BVP on $\caQ_\delta)$}
\label{thm:perturbed}
There are constants $C,M$ independent of $\delta\in ]0,1[$ 
and $\tau\in \{{\rm Re}\,\tau >M\}$
so that if 
$u\in H^2(\caQ_\delta)$ satisfies   the direct problem
\eqref{eq:PHBVP}  (resp. the transposed problem \eqref{eq:ADJBVP})
with $g_1=0$ and $g_2=0$ 
 then
\begin{equation}
\label{eq:h1est}
\begin{aligned}
|   \tau    |\, ({\rm Re}\,\tau) \,  & \| {u}   \|_{L^2(\caQ_\delta)}^2   
\ +\ 
\cr
 \frac{\Re\tau}{|\tau|}  &  \Big( 
  \|  |\beta|^{1/2}  {u}\|_{L^2(\partial\caQ_\delta)}^2
 +
\|\nabla_x {u}\|_{L^2(\caQ_\delta)}^2\Big)
 \ \le\
C\,
\Big|
\int_{\caQ_\delta}
f\ \ou\ dx
\Big|
\,.
\end{aligned}
\end{equation}
\end{theorem}

\begin{definition}
$\bullet$  Using the analytic continuation 
$H_\uOmega(X(\tau,x))$ from Part {\bf   iii}
of Proposition \ref{prop:analyticBC},
define $\Phi,\beta \in C^\infty(\{{\rm Re}\,\tau >M\} \times]0,1[\times \partial{\caQ_\delta})$ by
\begin{equation}
\label{eq:phibeta}
\begin{aligned}
&\Phi (\tau,x)\ :=\
\Pi(\tau,x)\,
  \Big(\sum_j
    \frac{\nu_j^2\tau^2}
   {(\tau+\sigma_j)^2}\Big)^{1/2},
   \cr
 &  \beta(\tau,\delta, x) 
\ :=\
 \tau + 2H_\uOmega\big(X(\tau,x)\big).
 \end{aligned}
      \end{equation}

$\bullet$ With ${\rm Re}\,\tau>M$ and $a(u,v)$ from  \eqref{eq:pform}, 
define  continuous bilinear forms $\caA(\tau, \cdot,\cdot ) :H^1({\caQ_\delta})\times H^1({\caQ_\delta}) \to \CC$ by
\begin{equation}
\label{eq:quadform}
\caA    (\tau, {u}    , {v})  \, :=\,
a({u}\,,\, {v})
\,+\,
\int_{{\caQ_\delta}} 
\tau^2\,\Pi(\tau,x)
\, u\cdot v\, dx
\, +\, 
\int_{\partial {\caQ_\delta}}
\Phi
\,
\beta\,
{u}
\,\cdot\,
 {v}
\,
d\Sigma.
\end{equation}
\end{definition}

The proof of Theorem \ref{thm:perturbed} relies on two estimates for  $\caA$.
The first is a lower
bound for $\caA(u,\ou)$ that holds for all $u\in H^1({\caQ_\delta})$.  The second is 
an upper bound
that 
relies  on the boundary conditions.
The dependence of $\caA$ on ${\caQ_\delta}$ and therefore $\delta$ is 
suppressed.  Similarly,
the dependence of $\caA$ on $\tau$ is usually
not indicated.
%
%

\begin{lemma}
\label{lem:VFpauli}
If $u\in H^2({\caQ_\delta})$ and  $v\in H^1({\caQ_\delta})$, define $f:= (\tau^2\Pi -p)u$.
Then,
\begin{equation}
\label{eq:Aidentity}
\begin{aligned}
\caA(\tau, {u}\,,{v}) \, -\,  \int_{{\caQ_\delta}} f\cdot v\,   dx
\, =\, 
\int_{\partial{\caQ_\delta}}
  \Phi(\tau,x)\ 
\big(V \, +\, \beta(\tau,\delta,x)\big) {u}
     \cdot
   {v} 
   \ d\Sigma\,.
\end{aligned}
\end{equation}
\end{lemma}

{\bf Proof.} 
  The differential operator appearing in the boundary term of
Green's formula
\eqref{eq:helmnormal}
is related to the operator $V(\tau,x,\partial)$ associated to the 
natural boundary condition for the stretched Pauli system by
\begin{equation*}
\label{eq:notV}
\begin{aligned}
\Pi(\tau,x)
   \sum_{j=1}^3
   \frac{\nu_j\,\tau^2}{   (\tau +\sigma_j)^2}\
   \partial_j
   \, &=\,
   \Pi(\tau,x)
  \Big( \sum_{j=1}^3
  \Big( \frac{\nu_j \,\tau^2}{   (\tau +\sigma_j)^2} \Big)^2\Big)^{1/2}  V(\tau,x,\partial)
  \cr
   \, &=\,
   \Phi(\tau,x)\,
    V(\tau,x,\partial).
  \end{aligned}
   \end{equation*}

Equation \eqref{eq:helmnormal} shows 
 that  
\begin{equation}
\label{eq:halfway}
\begin{aligned}
a & (   {u}\,,{v}) 
\ +\  
\int_{{\caQ_\delta}} 
\tau^2    \,\Pi(\tau,x)
\ u\cdot v\ dx
\ -\  \int_{{\caQ_\delta}} f\cdot v\, dx\ =\ 
\cr
 &
 \int_{\partial{\caQ_\delta}}
  \Pi(\tau,x)\ 
     \sum_{j=1}^3
   \frac{\nu_j \, \tau^2\,\partial_j  {u}}{   (\tau +\sigma_j)^2}\
   \cdot
   { {v} } \ d\Sigma
 \ =\  \int_{\partial{\caQ_\delta}}
  \Phi(\tau,x)\ 
V {u}
     \cdot
   {v} 
   \ d\Sigma
   \,.
    \end{aligned}
\end{equation}

Adding $\int_{\partial{\caQ_\delta}} \Phi\, \beta\, u\cdot v \,d\Sigma$  to both sides 
proves  \eqref{eq:Aidentity}.
\hfill
\qed
\vskip.2cm

\begin{example}
\label{ex:weakform}   If on $\partial{\caQ_\delta}$, $u$ satisfies
 $$
 \pi^-(\widetilde\nu)u=0,
 \quad
 {\rm and}
 \quad
 \Big(
V+ \tau+2H_\uOmega(X(\tau,\ux))
\Big)u=0,
$$
and $v$ satisfies $\pi^+(\widetilde\nu)^\dagger v=0$, then 
the boundary term in \eqref{eq:Aidentity}   vanishes.  
This yields a weak formulation, and a mixed finite element
approach to the boundary  value problem for $u$.
\end{example}

\subsubsection{Lower bound for $|\caA(u,\ou)|$}

\begin{proposition}
\label{prop:Alower}
There 
are constants $C,M>0$  independent of $\delta\in ]0,1[$
so that for any
$\tau\in \big\{\Re\tau\ge M\big\}$, and ${u}\in H^1(\caQ_\delta)$, 
\begin{equation}
\label{eq:coercivitypauli}
\begin{aligned}
 |\tau|(\Re\tau) &
 \| {u}\|_{L^2(\caQ_\delta)}^2
\  +\
 \cr
& \frac{\Re\tau}
 {|\tau|}\Big(
 \big\| |\beta|^{1/2}    \,u\big\|_{L^2(\partial\caQ_\delta)}^2
 +\
  \big\|
\nabla_x {u}\big\|_{L^2(\caQ_\delta)}^2  \Big)  
 \lesssim\
 \big|   \caA  ( {u}, \ou)  \big|
 .
 \end{aligned}
\end{equation}
\end{proposition}

\begin{remark}     In \eqref{eq:zinfty}, we show that $H_\uOmega = H_{\caQ_\delta} + O(1/\tau)$.
Since $\beta=\tau + 2\,H_{\uOmega}(\tau,x) $ 
it follows that there is an $M$ independent of $\delta$ to that
for $\Re \tau >M$ 
$$
 |\tau| +H_{\caQ_\delta}(x)
\ \le\
 |\beta(\tau,\delta ,x)|
 \ \le\
  |\tau| +3H_{\caQ_\delta}(x) \,.
 $$
\end{remark}

{\bf Proof.}   {\bf Step 1. $ \caA_0$ and its real and imaginary parts.}   
Denote by $\caA_0$ the 
form that one would have if $\sigma_j=0$ for all $j$,
\begin{align*}
\caA_0(\tau, {u},  v)
\ &:=\
 \int_{\caQ_\delta}
 \tau^2\ 
 u\cdot v\ dx\ +\
 \int_{\partial\caQ_\delta}
 \beta  
\,
u\cdot v
\
d\Sigma
 \ +\
  \int_{\caQ_\delta}
 \nabla_x{u}
 \cdot
 \nabla_x v \ dx\,,
 \cr
\caA_0(\tau, {u},  \ou)
\ &:=\
 \int_{\caQ_\delta}
 \tau^2\ 
 | {u}|^2\ dx\ +\
 \int_{\partial\caQ_\delta}
  \beta
\
|{u}|^2
\
d\Sigma
 \ +\
  \int_{\caQ_\delta}
 \big|\nabla_x{u}\big|^2\ dx\,.
\end{align*}
 The real part  of $\caA_0$
 is 
 \begin{equation}
 \label{eq:recoer}
 \begin{aligned}
 \Re \caA_0( {u}, \ou)
 \ =\
  \Big((\Re\tau & )^2  -({\Im}\,  \tau )^2  \Big)
 \big\| {u}\big\|_{L^2(\caQ_\delta)}^2
\cr
\ +\ & \,
\big\|
(\Re \beta)^{1/2}
{u}
\|^2_{L^2(\partial\caQ_\delta)}
+
\big\|\nabla_x {u}\big\|_{L^2(\caQ_\delta)}^2\,.
 \end{aligned}
 \end{equation}

 Use $\,{\Im}\,\tau^2=2\ ({\Im}\,\tau) ({\Re}\,\tau)$ 
to find,
\begin{equation*}
\begin{aligned}
 {\Im} 
 \int_{\caQ_\delta}
 \tau^2\ 
 |{u}|^2\ dx
\ &=\
({\Im}\, \tau)
\int_{\caQ_\delta}
 2\, {\Re}\, \tau
 \
 | {u}|^2\, dx \,,
 \cr
\Im
\int_{\partial\caQ_\delta}
\beta
\
| {u}|^2
\
d\Sigma
\ &=\
(\Im \tau) 
\int_{\partial\caQ_\delta}
\
| {u}|^2
\
d\sigma
\,.
\end{aligned}
\end{equation*}

Combining shows that for $0\ne \Im\tau$,
\begin{equation}
\label{eq:imcoer}
\begin{aligned}
\frac{
\Im \caA_0({u},  \ou)  }
{\Im\tau}
\ =
\
2\, (\Re \tau)\,
 \big\| {u}\big\|_{L^2(\caQ_\delta)}^2
 +
  \big\|  {u}\big\|_{L^2(\partial\caQ_\delta)}^2
\, .
\end{aligned}
\end{equation}

 {\bf Step 2.  Proof for $\caA_0$.}
$\bullet$   The bound
 \eqref{eq:coercivitypauli} 
  is proved by combining \eqref{eq:recoer} and
 \eqref{eq:imcoer}.  Care is needed 
 where the terms on the right of \eqref{eq:recoer}
 do not  have the same sign.
\JR{
On the set
 $\big\{|\Im \tau| <  \Re \tau/2\big\}$,  \eqref{eq:recoer}
 implies  
 \eqref{eq:coercivitypauli}.
 }

$\bullet $   It suffices to consider  the complementary set $\{|\Im \tau|\ge \Re \tau/2\}$.
 In that parameter range \eqref{eq:imcoer} implies
 \begin{equation}
\label{eq:imag}
 (\Re \tau)\,
 \| {u}\|_{L^2(\caQ_\delta)}^2
\ +\
  \| {u}\|_{L^2(\partial\caQ_\delta)}^2
\   \lesssim\ 
  \frac{|\Im  \caA_0(u,\ou)|}
  {| \tau|}
  \ .   
  \end{equation}
     Multiplying by $| \tau|^2/\Re\tau$ yields
  \begin{equation}
  \label{eq:other}
|\tau|^2\, 
 \| {u}\|_{L^2(\caQ_\delta)}^2
\ +\
 \frac{ |\tau|^2}{\Re\tau}\,
    \| {u}\|_{L^2(\partial\caQ_\delta)}^2
 \    \lesssim\
 \frac{ |\tau|}{\Re\tau}\,   |\Im  \caA_0( u, \ou)|
  \,.
     \end{equation}
  Therefore,
$$
  \Big|(\Re\tau)^2  -({\Im}\,  \tau )^2  \Big|\
 \big\| {u}\big\|_{L^2(\caQ_\delta)}^2
 \ \le\
 |\tau|^2\, 
 \| {u}\|_{L^2(\caQ_\delta)}^2
 \ \lesssim\
  \frac{ |\tau|}{\Re\tau}\, 
\big|
 \caA_0(  {u},    \ou)
\big|
\,.
$$
   Using this in \eqref{eq:recoer} yields
  for $|\tau|>M_1$,
    \begin{equation}
 \label{eq:domi}
  \big\|
\nabla_x {u}\big\|_{L^2(\caQ_\delta)}^2
 +
 \big\| (\Re\beta|)^{1/2}    \,u\big\|_{L^2(\partial\caQ_\delta)}^2
\ \lesssim\
 \frac{ |\tau|}{\Re\tau}\,  \big| \caA_0( {u}, \ou)  \big|
 \,.
 \end{equation}
 Adding \eqref{eq:other} and \eqref{eq:domi} yields
      \begin{align*}
 |\tau|^2\, 
 \| {u}\|_{L^2(\caQ_\delta)}^2
\ +
\
  \big\|
\nabla_x {u}\big\|_{L^2(\caQ_\delta)}^2
\ &+\
 \big\| (\Re \beta)^{1/2}    \,u\big\|_{L^2(\partial\caQ_\delta)}^2
 \ +\
 \cr
 \ &+\ 
 \frac{ |\tau|^2}{\Re\tau}\,
    \| {u}\|_{L^2(\partial\caQ_\delta)}^2
\ \lesssim\
 \frac{ |\tau|}{\Re\tau}\,  \big| \caA_0( {u}, \ou)  \big|
 \,.
 \end{align*}
 Multiply by $(\Re\tau)/|\tau|$ 
 and use $|\beta|\le (\Re \beta) + |\tau|$ 
 to find
\JR{the desired estimate},
 \begin{equation}
\label{eq:domi2}
\begin{aligned}
 |\tau|(\Re\tau) 
 \| {u}\|_{L^2(\caQ_\delta)}^2
\,  +\,
 \frac{\Re\tau}
 {|\tau|}\Big(
 \big\| |\beta|^{1/2}    \,u\big\|_{L^2(\partial\caQ_\delta)}^2
 +\
  \big\| 
\nabla_x &  {u}\big\|_{L^2(\caQ_\delta)}^2  \Big)  
\cr
\  \lesssim\  &
 \big|   \caA_0  ( {u}, \ou)  \big|
 .
 \end{aligned}
\end{equation}

  \vskip.1cm
   
   {\bf Step 3.  Perturbation argument.}   For 
   $\tau\ne 0$,
   $
   \tau + \sigma_j(x_j)  =
   \tau(
   1+ 
   \sigma_j/\tau)
   )
   $.
   Write
   \begin{align*}
   a(u,\ou) &- a_0(u, \ou) 
   =
   \int_{\caQ_\delta}
   \Big(
   \frac{
   (\tau + \sigma_{j+1})
    (\tau + \sigma_{j+2})}
    {\tau( \tau + \sigma_{j})}
   -1
   \Big)
   |\partial_j u|^2\ dx
   \cr
   &+
   \tau^2
   \int_{{\caQ_\delta}}
   \big(\Pi(\tau,x)-1\big)
   |u|^2\ dx
   +
   \int_{\partial{\caQ_\delta}}
   \big(
   \Phi(\tau,x)-1\big)
   |\beta|\,|u|^2
   \ d\Sigma\,.
\end{align*}
Since $|\Pi-I|+|\Phi-I|=O(1/\tau)$,
 this yields
         \begin{equation}
   \label{eq:isabella}
   \begin{aligned}
   \big|  
   \caA & (  {u}, \ou) -\caA_0( {u}, \ou)
   \big|
 \  \lesssim\
   |\tau|   \| {u}\|_{L^2(\caQ_\delta)}^2+
   \cr &
 +
 \frac1{|\tau|}
  \||\beta|^{1/2} {u}\|_{L^2(\partial\caQ_\delta)}^2
+
\frac1{|\tau|}
\|\nabla_x {u}\|_{L^2(\caQ_\delta)}^2
  \, \lesssim\,
  \frac1{\Re\tau}\
  \big|\caA_0( {u}\,,\,  \ou) 
  \big|,
   \end{aligned}
   \end{equation} 
   where inequality \eqref{eq:domi2}
   for $\caA_0$ is used in the last inequality.
   The triangle inequality and     estimate   \eqref{eq:isabella}   imply
   \begin{equation*}
   \begin{aligned}
  \big| \caA( {u}\,,\, \ou)\big|
    \, \ge\, 
    \caA_0({u}\,,\,\ou)
    \,  -\,
    \big|
    \caA( {u}\,,\, \ou)
    -
    \caA_0({u}\,,\, \ou)
    \big|
    \, \ge\, 
     \Big(
  1 - \frac{c}{\Re\tau}
   \Big)\
  \big|
    \caA_0 ( {u}\,,\, \ou)\big|.
     \end{aligned}
   \end{equation*}

     For $\Re \tau>2c$ this
     yields  \eqref{eq:coercivitypauli} completing the proof of Proposition
   \ref{prop:Alower}.
   \hfill
\qed
\vskip.3cm

\subsubsection{Upper bound for $|\caA(u,\ou)|$,  proof of Theorem \ref{thm:perturbed}}
\label{sec:UP}

\begin{proposition}
\label{prop:Aupper}
If $u\in H^2({\caQ_\delta})$ is a solution of the Helmholtz boundary value problem
\eqref{eq:PHBVP} (resp.  the transposed problem \eqref{eq:ADJBVP}) with $g_1=0$ and $g_2=0$, 
then with constant independent of $\delta\in ]0,1[$ and $|\tau|>1$,
\begin{align*}
\big|\caA(u,\overline u)  & \big|
\,  \lesssim\, 
\Big|
\int_{\caQ_\delta} f\ \ou\ dx\Big|
\, +\,
  \frac1{ |\tau|}\,
 \Big(
  \, \big\||\beta|^{1/2} \,u\big\|^2_{L^2(\partial{\caQ_\delta})}
+
 \| u\|^2_{H^1({\caQ_\delta})}
 \Big)  \,  .
\end{align*}
\end{proposition}

{\bf Proof of Proposition \ref{prop:Aupper}.}  
To treat  \eqref{eq:PHBVP},
 write
$$
\big(V +\beta\big) u
=
(\pi^+(\widetilde \nu) 
+
\pi^-(\widetilde \nu) )
\big(V+\beta\big) u
\ =\
\pi^-(\widetilde \nu )
\big(V+\beta\big) u.
$$
For the transposed  boundary value problem \eqref{eq:ADJBVP} write
$$
\big(V +\beta\big) u
\ =\ 
(\pi^+(\widetilde \nu)^\dagger 
+
\pi^-(\widetilde \nu)^\dagger )
\big(V +\beta\big) u
\ =\
\pi^-(\widetilde \nu )^\dagger
\big(V+\beta\big) u.
$$

Continuing the computation for \eqref{eq:PHBVP},
Lemma \ref{lem:VFpauli} yields for $u,v\in H^1({\caQ_\delta})$,
\begin{equation*}
\label{eq:halfgreen3}
\begin{aligned}
\caA(u,v) 
\, =\, 
\int_{\caQ_\delta}
f \,\cdot\,      v\, dx 
\, -\, 
\int_{\partial{\caQ_\delta}} \Phi\,  \pi^-(\widetilde \nu) \, \big( V(\tau,x,\partial)+\beta\big) u\,\cdot\, v\ d\Sigma.
\end{aligned}
\end{equation*}
With $v=\ou$ this is
\begin{equation}
\label{eq:halfgreen4}
\begin{aligned}
\caA(u,\ou) 
\, =\, 
\int_{\caQ_\delta}
f \,\cdot\,      \ou\, dx 
\, -\, 
\int_{\partial{\caQ_\delta}} \Phi\,  \pi^-(\widetilde \nu) \, \big( V(\tau,x,\partial)+\beta\big) u\,\cdot\, \ou\ d\Sigma.
\end{aligned}
\end{equation}

The difficult step  is to derive an upper bound for 
$$
\int_{\partial{\caQ_\delta}} \Phi\,  \pi^-(\widetilde \nu)\big( V(\tau,x,\partial)+\beta\big) u\,\cdot\, \overline u\ d\Sigma.
$$
The boundary condition  $\pi^-(\widetilde \nu)u=0$ implies $\pi^+(\widetilde \nu)u=u$ so,
$$
\int_{\partial{\caQ_\delta}} \Phi\,  \pi^-(\widetilde \nu)\big( V(\tau,x,\partial)+\beta\big) u\,\cdot\, \overline u\, d\Sigma
=
\int_{\partial{\caQ_\delta}} \Phi\,  \pi^-(\widetilde \nu)\big( V(\tau,x,\partial)+\beta\big) u\,\cdot\, \overline{\pi^+(\widetilde\nu) u}\, d\Sigma.
$$
Write
$$
\overline{\pi^+(\widetilde\nu) u}
\,=\,
\overline{\pi^+(\widetilde\nu)}  \, \overline u
\,=\,
\pi^+(\widetilde\nu)^\dagger  \,\overline u
\,+\,
\Big(\overline{\pi^+(\widetilde\nu) }
-\pi^+(\widetilde\nu)^\dagger \Big)
\overline u
\,.
$$
When this is inserted the $(\pi^+(\widetilde \nu))^\dagger \overline u$ term yields zero.  Therefore
\begin{equation}
\begin{aligned}
\label{eq:bad}
\int_{\partial{\caQ_\delta}}  & \Phi\, \pi^-(   \widetilde \nu)\big(  V(\tau,x,\partial)    +\beta\big) u\,\cdot\, \overline u\ d\Sigma
\cr
&=
\int_{\partial{\caQ_\delta}} \Phi\,  \pi^-(\widetilde \nu)\big( V(\tau,x,\partial)+\beta\big) u\,\cdot\, 
\Big(\overline{\pi^+(\widetilde\nu) }
-\pi^+(\widetilde\nu)^\dagger \Big)
\overline u
\ d\Sigma
\cr
&=
\int_{\partial{\caQ_\delta}} \Phi\, \big( V(\tau,x,\partial)+\beta\big) u\,\cdot\, 
 \pi^-(\widetilde \nu)^\dagger
\Big(\overline{\pi^+(\widetilde\nu) }
-\pi^+(\widetilde\nu)^\dagger \Big)
\overline u
\ d\Sigma
\cr
&=
\int_{\partial{\caQ_\delta}} \Phi\, \big( V(\tau,x,\partial)+\beta\big) u\,\cdot\, 
w
\ d\Sigma
\end{aligned}
\end{equation}
with
$$
w:=  \pi^-(\widetilde \nu)^\dagger
\Big(\overline{\pi^+(\widetilde\nu) }
-\pi^+(\widetilde\nu)^\dagger \Big)
\overline u.
$$

For the transposed  problem 
the difficult boundary term is 
$$
\int_{\partial{\caQ_\delta}}   
\Phi\, \pi^-(   \widetilde \nu)^\dagger\big(  V(\tau,x,\partial)    +\beta\big) u\,\cdot\, \overline u\ 
d\Sigma=
\int_{\partial{\caQ_\delta}} \Phi\, \big( V(\tau,x,\partial)+\beta\big) u\,\cdot\, 
w
\ d\Sigma
$$
with
$$
w:=  \pi^-(\widetilde \nu)
\Big(\overline{\pi^+(\widetilde\nu)^\dagger }
-\pi^+(\widetilde\nu) \Big)
\overline u.
$$

The estimates in the two cases are  virtually identical.  The  details
are presented only for the direct problem.
For the direct problem define $\sfm\in C^\infty(\{{\rm Re}\,\tau >M\}\times \partial {\caQ_\delta}$ by 
\begin{equation}
\label{eq:defm}
\sfm(\tau,x) := \tau\ \pi^-(\widetilde \nu)^\dagger\
\big(\overline{\pi^+(\widetilde\nu) }
-\pi^+(\widetilde\nu)^\dagger \big),
\qquad
{\rm so},
\quad
w = \frac1{\tau}\, \sfm \, \ou\,.
\end{equation}
\JR{ Equation \eqref{eq:defm} shows 
 that
  \eqref{eq:bad}  is equal to }
\begin{equation}
\label{eq:bad2}
\frac1{\tau}\Big(
\int_{\partial{\caQ_\delta}} \Phi Vu\cdot \sfm \ou\ dx
\, +\,
\Phi\ \beta\ u\cdot \sfm \, \ou\ d\Sigma
\Big).
\end{equation}
The next lemma gathers estimates for $V$ and $\sfm$.

\begin{lemma}
\label{lem:alexis}
There are constants $C,M$ so that for all
${\rm Re}\,\tau>M$, and, $0<\delta<1$,
the following hold.

${\bf i.} \ \ \ {\rm supp}\,   \sfm \ \subset\
\big\{ x\in \partial{\caQ_\delta}\ :\ {\rm dist}(x,\caS)<\delta\big\} $.

${\bf ii.}\ \ \| \sfm(\tau,x)\|_{L^\infty(\partial{\caQ_\delta})}\     \le\
C$.

${\bf iii.}\ \big\| \nabla_x \sfm(\tau,x)\big\|_{L^\infty(\partial{\caQ_\delta})}\ \le\
C\,
|\beta|
\,.
$

{\bf iv.}    For all $u\in H^{1/2}(\partial{\caQ_\delta})$,
$$
\|\sfm\, u\|_{H^{1/2}(\partial {\caQ_\delta})}
\ \lesssim\
\| \,|\beta(\tau,x)|^{1/2}\,
u
\|_{L^{2}(\partial {\caQ_\delta})}
+
 \|u\|_{H^{1/2}(\partial {\caQ_\delta})}.
 $$

\vskip.1cm
{\bf v.} \  For all $u\in H^1({\caQ_\delta})$,
$
\|Vu\|_{H^{-1/2}(\partial{\caQ_\delta})}
\ \le \ 
C\,
 \| u\|_{H^{1}({\caQ_\delta})}.
$

\end{lemma}

{\bf Proof of Lemma.}  {\bf i.}  
   For most points $x\in \partial{\caQ_\delta}$,  one has
 $x\in G_j$ for some $j$, $\nu=\pm\bfe_j$, and $A(\nu)=\pm A_j$ 
 The spectral representation is 
$$
A(\nu) = \pi^+(\nu) -\pi^-(\nu),
\quad
\pi^{\pm}(\nu) = \pi^{\pm}(\nu)^* = \pi^{\pm}(\nu)^2,
\quad
\pi^{\pm}(\nu) \pi^{\mp}(\nu)=0. 
$$
  These imply the spectral representations
  $$
   A(\nu)^\dagger \ =\  \pi^+(\nu)^\dagger\ -\ \pi^-(\nu)^\dagger,
   \qquad
   {\rm and},
   \qquad
   \overline{A(\nu)}\ =\  \overline{\pi^+(\nu)}\  -\  \overline{\pi^-(\nu)}
   \,.
   $$
   
   For $j\in\{1,2, 4,5\}$, $A(\nu)$ is real and hermitian symmetric, 
   $A(\nu) = A(\nu)^\dagger=\overline{A(\nu)}$.  
 Comparing the spectral representations 
 yields
 $ \pi^\pm (\nu)  = 
 \pi^\pm(\nu)^\dagger
  =
 \overline{\pi^\pm(\nu)}$.  Since $\widetilde \nu$ is a scalar multiple of $\nu$ this
 yields 
 \begin{equation}
 \label{eq:longpi}
 \pi^\pm (\nu) \ =\ 
 \pi^\pm(\nu)^\dagger
 \ =\
 \overline{\pi^\pm(\nu)}
 \ = \
 \pi^\pm(\widetilde \nu)
 \ =\
 \pi^\pm(\widetilde \nu)^\dagger
 \ =\
 \overline{  \pi^\pm(\widetilde \nu)}
 \,.
 \end{equation}
It follows that $\sfm=0$ at such points.  

For $j\in \{3,6\}$ and $x\in G_j$, $A(\nu)^\dagger = \overline {A(\nu)} = -A(\nu)$.  Comparing
the spectral 
representations as above implies that 
 \begin{equation}
 \label{eq:longpi2}
 \pi^{\pm}(\nu) \ =\   \pi^{\mp}(\nu)^\dagger \ =\  \overline{\pi^{\mp}(\nu)}
 \ =\
  \pi^{\pm}(\widetilde\nu) \ =\   \pi^{\mp}(\widetilde\nu)^\dagger \ =\  \overline{\pi^{\mp}(\widetilde \nu)}.
 \end{equation}
Therefore 
$\sfm=0$ at these flat parts of the boundary too.  

These results for all $G_j$ show that $\sfm$ is supported
on the rounded edges of $\partial{\caQ_\delta}$ proving {\bf i.}

\vskip.1cm
$\bf ii.$  Compute
$$
\frac{\tau}{\tau + \sigma_j}
\,=\,
\frac{1}{1 + \sigma_j/\tau}
\,=\,
1 -\frac{\sigma_j}
{\tau}
+
\Big(\frac{\sigma_j}
{\tau }\Big)^2 -\ \cdots\,.
$$
It follows that as $|\tau|\to \infty$,
$$
\widetilde \nu - \nu  \,=\, O(1/|\tau|),
\quad
{\rm so},
\quad
\pi^+(\widetilde\nu) -\pi^+(\nu) 
\,=\, O(1/|\tau|).
$$
To estimate the size of $\sfm$ write
$$
\overline{\pi^+(\widetilde\nu) }
-\pi^+(\widetilde\nu)^\dagger 
\,=\,
\Big(
\overline{\pi^+(\widetilde\nu) }
-
\overline{\pi^+(\nu)}
\Big)
\,+\,
\Big(
\overline{\pi^+(\nu) }
-
\pi^+(\widetilde\nu)^\dagger
\Big)\,.
$$
The first summand is $O(1/\tau)$.  
Equations \eqref{eq:longpi} and
\eqref{eq:longpi2}
imply  that the second is
also
$O(1/|\tau|)
$.
If follows that  $\sf m$ is bounded uniformly in $\tau,\delta$, proving
{\bf ii.}

\vskip.1cm
$\bf iii.$   Use the notations from 
Proposition \ref{prop:analyticBC}.  Then $\tau\mapsto \nu(\tau,\cdot )$ is analytic
in $|\tau|>R$ with  values in $C^\infty(\partial{\caQ_\delta})$.

  Expand the stretchings  in $z=1/\tau$ about $z=0$.  The 
  transformation satisfies
   \begin{equation}
  \label{eq:Xj2}
  \frac{\partial X_j(\tau,x_j)}{\partial x_j}   \ =\ 
  \frac{\tau + \sigma_j(x_j)}
  {\tau}
  \ =\ 
  {1 + z\,\sigma_j(x_j)}
  \,,
  \qquad
  X_j(\tau,0)\,=\,0\,.
  \end{equation}
  Thus $X$ is analytic on a neighborhood of $z=0$ with 
  $X(0,x) =x$.  The derivative with respect to $x$ satisfies
  $D_xX = I + O(z)$.   It follows that 
  $$
  \nu(\tau, x) = \nu(\infty,x) + O(1/\tau),
  \quad
  {\rm and},
  \quad
  \nabla_x\nu(\tau, x) 
  =
   \nabla_x\nu(\infty, x) +
   O(z).
  $$ 
At $\tau=\infty$
the $\nabla_x\nu$ restricted to the tangent space is the Weingarten map of $\partial{\caQ_\delta}$
from Definition \ref{def:Weingarten}.
At $\tau=\infty$, the eigenvalues are nonnegative.   Therefore
\begin{equation}
\label{eq:zinfty}
\begin{aligned}
&H_\uOmega(\tau,x)\ =\  H_{\caQ_\delta}(x) + O(1/\tau),
\cr
&|  \nabla_x\nu(\infty,x) |
\ \lesssim\ 
\max\,\{\kappa_1,\kappa_2\} 
\ \le\
 2H_{\caQ_\delta}(x),
 \cr
 &|  \nabla_x\nu(\tau,x) |
\,\lesssim\,
  |H_\uOmega(\tau,x) | + |\tau|^{-1} 
  \,\lesssim\,
| \beta(\tau,\delta,x)|  \,.
 \end{aligned}
\end{equation}
Since
$
|\nabla_x\sfm|
\lesssim
 |\nabla_x\nu|
$
this proves
{\bf iii.}

\vskip.1cm
{\bf iv.}  Estimates {\bf ii, iii}  imply that with constants independent
of $\tau,\delta$ and all $u$,
\begin{equation}
\begin{aligned}
\label{eq:2m}
\|\sfm\, u\|_{L^2(\partial {\caQ_\delta})}
&\lesssim \, \| u\|_{L^2(\partial {\caQ_\delta})},
\cr
\|\sfm\, u\|_{H^1(\partial {\caQ_\delta})}
&\lesssim
\|\,
| \beta  | \,
  u\|_{L^2(\partial {\caQ_\delta})} 
  \ +\ 
  \|  u\|_{H^1(\partial {\caQ_\delta})} \,.
\end{aligned}
\end{equation}
To prove  the second, apply the product rule  with vector fields $\partial$
that are tangent to the boundary to find
$\partial(\sfm u) = \sfm \partial u + (\partial\sfm)u$.
Therefore
$$
\|\partial(\sfm u)\|_{L^2(\partial{\caQ_\delta})}
\ \le\ 
\|\sfm\|_{L^\infty(\partial{\caQ_\delta})}
\| \partial u\|_{L^2(\partial{\caQ_\delta})}
\ +\ 
\|(\partial \sfm)\, u\|_{L^2(\partial{\caQ_\delta})}
\,.
$$
Using {\bf iii} in the second summand proves \eqref{eq:2m}

Denote by $\Delta_{\partial{\caQ_\delta}}$ the Laplace-Betrami
operator of $\partial{\caQ_\delta}$.
The estimates  \eqref{eq:2m} are the cases $\theta=0,1$ of 
\begin{align*}
\|\sfm\, u\|_{H^{\theta}(\partial {\caQ_\delta})}
 \lesssim
  \|
  (
  | \beta(\tau,x) | + |\Delta_{\caQ_\delta}|^{1/2}   )^\theta
  u\|_{L^2(\partial {\caQ_\delta})} 
  .
  \end{align*}
  Interpolation implies the estimate for 
$0\le \theta\le 1$.   Use  the case $\theta=1/2$.
For self adjoint $B_j\ge0$ with  $B_1$ bounded and  $u\in \caD(B_2)$,
\begin{align*}
\big\|
\sqrt{B_1+B_2}\, u\big\|^2
\ &=\
\big( 
\sqrt{B_1+B_2}\, u,
\sqrt{B_1+B_2}\, u
\big)
\ =\ ((B_1+B_2)u,u)
\cr
\ &=\ 
(B_1u,u) + 
(B_2u,u)
\ =\
\big\|
\sqrt{B_1}\, u\big\|^2
 +
\big\|
\sqrt{B_2}\, u\big\|^2
\,.
\end{align*}
With $B_1=|\beta(\tau,x)|$ and $B_2= |\Delta_{\caQ_\delta}|^{1/2}$ this yields
$$
 \|
  (
  | \beta(\tau,x) | +  |\Delta_{\caQ_\delta}|^{1/2} )^{1/2}
  u\|_{L^2(\partial {\caQ_\delta})} ^2
  \, =\,
   \|
  | \beta(\tau,x) |^{1/2}
  u\|_{L^2(\partial {\caQ_\delta})} ^2
  \ +\
  \| |\Delta_{\caQ_\delta}|^{1/4}   u\|_{L^2(\partial {\caQ_\delta})}^2.
  $$
Using  this  in the $\theta=1/2$ estimate proves
{\bf iv.}

\vskip.1cm
{\bf v.}
With constants independent of $\delta,\tau$ with $|\tau|>R$,
one has for all $u\in H^1({\caQ_\delta})$,
$$
\int_{\caQ_\delta}
\big| \nabla_xu\big|^2\ dx
\ \le \ 
C\Big( -{\rm Re}
\int_{\caQ_\delta}
p(\tau,x,\partial)\,u \cdot \ou\ dx
\Big).
$$
It follows that 
for $|\tau|>R$ and $0<\delta<1$,
the operator
$1-p(\tau,x,\partial)$ is an isomorphism of $H^1({\caQ_\delta})$ to $H^{-1}_0({\caQ_\delta})$,
and with 
constants independent of $\tau, \delta$,
$$
\|u\|_{H^{1}({\caQ_\delta})}
\ \lesssim\
\|(1-p)u\|_{H^{-1}_0({\caQ_\delta})}
\ \lesssim\
\|u\|_{H^{1}({\caQ_\delta})}\,.
$$
Therefore,
\begin{equation}
\label{eq:pu}
\begin{aligned}
\|pu\|_{H^{-1}_0({\caQ_\delta})}
\ &\le \
\|(1-p)u\|_{H^{-1}_0({\caQ_\delta})}
+
\|u\|_{H^{-1}_0({\caQ_\delta})}
\cr
\ &\lesssim\
\|u\|_{H^{1}({\caQ_\delta})} 
+
\|u\|_{H^{-1}_0({\caQ_\delta})}
\ \lesssim\
\|u\|_{H^{1}({\caQ_\delta})} \, .
\end{aligned}
\end{equation}

 Using \eqref{eq:defV}, \eqref{eq:helmnormal}, and 
 \eqref{eq:phibeta} shows that
 for all $u,v\in H^1({\caQ_\delta})$
$$
a(u,v) \ -\ 
\int_{\caQ_\delta} \big(p(\tau,x,\partial \big)\,u)\cdot v\ dx
\ =\ 
\int_{\partial{\caQ_\delta}}
\Phi(\tau,x)\ 
Vu \cdot v \ d\Sigma\,.
$$
For $\phi\in H^{1/2}(\partial{\caQ_\delta})$
choose $v\in H^1({\caQ_\delta})$ with 
$\|v\|_{H^1({\caQ_\delta})} \lesssim \|\phi\|_{H^{1/2}(\partial {\caQ_\delta})}$
to find,
\begin{align*}
\Big|\int_{\partial{\caQ_\delta}}
\Phi(\tau,x)\ 
Vu \cdot \phi\ d\Sigma
\Big|\ &=\
\Big|a(u,v) \ -\  
\int_{\caQ_\delta} (pu)\cdot v\ dx\Big|
\cr
\ &\lesssim\
\|\nabla u\|_{L^2({\caQ_\delta})}\,
\|\nabla v\|_{L^2({\caQ_\delta})}
\ +\ 
\|pu\|_{H^{-1}_0({\caQ_\delta})}
\|v\|_{H^{1}({\caQ_\delta})}
\cr
\ &\lesssim\
\big(
\|\nabla u\|_{L^2({\caQ_\delta})}
\,+\,
\|pu\|_{H^{-1}_0({\caQ_\delta})}
\big)\, \|\phi\|_{H^{1/2}(\partial{\caQ_\delta})}\,.
\end{align*}

Using this in 
the upper bound for  $|\int\Phi Vu\cdot\phi\,d\Sigma|$,
 shows that 
$$
\Big|\int_{\partial{\caQ_\delta}}
\Phi(\tau,x)\ 
Vu \cdot \phi\ d\Sigma
\Big|
\ \lesssim\
\|u\|_{H^{1}({\caQ_\delta})} \
\|\phi\|_{H^{1/2}(\partial{\caQ_\delta})}.
$$
Since $\Phi$ and $1/\Phi$ as well as their derivatives are
uniformly bounded, this proves {\bf v.}
\hfill
\qed
\vskip.2cm

{\bf End of proof of Proposition \ref{prop:Aupper}.}   The second term on the right
in \eqref{eq:bad2}  is estimated as
\begin{equation}
\label{eq:second}
\Big|
\int_{\partial{\caQ_\delta}} \Phi\ \beta\ u\cdot \sfm \, \ou\ d\Sigma
\Big|
\ \lesssim\  
\int_{\partial{\caQ_\delta}}  |\beta|\ |u|^2\ d\Sigma
\ =\
\| |\beta|^{1/2} u\|^2_{L^2(\partial{\caQ_\delta})}
\,.
\end{equation}

The first summand is estimated as
\begin{align}
\label{eq:start}
\Big|
\int_{\partial{\caQ_\delta}} \Phi Vu\cdot \sfm \ou\ d\Sigma
\Big|
\ &\lesssim\
\|Vu\|_{H^{-1/2}(\partial{\caQ_\delta})}
\
\|\sfm \ou\|_{H^{1/2}(\partial{\caQ_\delta})}\,.
\end{align}
For $\|\sfm \, u\|_{H^{1/2}(\partial{\caQ_\delta})} $ 
use Part {\bf iv} of the lemma in \eqref{eq:start}  
to find,
\begin{equation}
\label{eq:prec}
\begin{aligned}
\Big|
\int_{\partial{\caQ_\delta}}  & \Phi   Vu\cdot \sfm \ou\ d\Sigma
\Big|
\ \lesssim\ 
\cr
&\big(
\| p\, u \|_{H^{-1}_0({\caQ_\delta})}
+
\|\nabla u\|_{L^2({\caQ_\delta})}
\big)
\big(
 \| \,|\beta|^{1/2}u\|_{L^2(\partial{\caQ_\delta})}
+
 \| u\|_{H^{1/2}(\partial{\caQ_\delta})}\big).
\end{aligned}
\end{equation}

Use this, \eqref{eq:pu}, 
and, $ \| u\|_{H^{1/2}(\partial{\caQ_\delta})}\lesssim \|u\|_{H^{1}({\caQ_\delta})}$
in \eqref{eq:prec} to find,
\begin{align}
\label{eq:third}
\Big|
\int_{\partial{\caQ_\delta}}  & \Phi   Vu\cdot \sfm \ou\ d\Sigma
\Big|
\ \lesssim\
 \| \,|\beta|^{1/2}u\|_{L^2(\partial{\caQ_\delta})}^2
+
 \| u\|_{H^{1}({\caQ_\delta})}^2\,.
\end{align}

Adding the estimates
\eqref{eq:second},
and 
\eqref{eq:third}
 for the  terms on the right 
of \eqref{eq:halfgreen4}  proves Proposition \ref{prop:Aupper}.
\hfill
\qed
\vskip.2cm

{\bf Proof of Theorem \ref{thm:perturbed}.}  Combine the lower and upper bounds
for $|\caA(u,\ou)|$ from  Propositions \ref{prop:Alower} and \ref{prop:Aupper}
 to find,  
\begin{align*}
|   \tau    |\, ({\rm Re}\,\tau) \,   \| {u}   \|_{L^2(\caQ_\delta)}^2   
& + 
 \frac{\Re\tau}{|\tau|}\Big( 
  \|  |\beta|^{1/2}  {u}\|_{L^2(\partial\caQ_\delta)}^2
 +
\|\nabla_x {u}\|_{L^2(\caQ_\delta)}^2\Big)\cr
 &
 \le
C\,
\Big|
\int_{\caQ_\delta}
f\ \ou\ dx
\Big|
+
  \frac{C}{ |\tau|}\,
 \Big(
  \|   |\beta|^{1/2}  u\|^2_{L^2(\partial\caQ_\delta)}
+
 \|u\|^2_{L^2(\caQ_\delta)}
 \Big).
\end{align*}
Choose $M=2C$.   Then for ${\rm Re}\,\tau>M$,
the second summand on the 
right can be absorbed in the left  hand side yielding
\eqref{eq:h1est}.   This completes the proof of 
Theorem  \ref{thm:perturbed}.
\hfill
\qed
\vskip.2cm

\subsubsection{Proof of Theorem \ref{thm:taufixed}}
\label{sec:deltafixed}

\JR{This section carries out the strategy outlined} after the  statement
of Theorem \ref{thm:taufixed}.

\vskip.1cm
{\bf Proof that the map $u\mapsto (f,g_1,g_2)$ has trivial 
kernel.}    If $u\in C^\infty (\overline{\caQ_\delta})$ is in the kernel,
it follows that 
$u\in H^2(\caQ_\delta)$ 
and satisfies  the homogeneous boundary  value problem with sources 
$f,g_1,g_2$ equal to zero.
Theorem  \ref{thm:perturbed} implies 
that $u=0$.

\vskip.1cm
{\bf Proof that the  annihilator of  the range,
is $\{0\}$.}  
$\bullet$  \JR{Use} the following Green's identity for 
 $u ,v \in H^2({\caQ_\delta})$, 
\begin{align}
\label{eq:dotGreen}
& \int_{{\caQ_\delta}} \big(  
 \tau^2\Pi(\tau,x)   -p(\tau  ,x,\partial)  \big) u  \cdot v\ dx 
  -
 u  \cdot 
  \big(  
 \tau^2\Pi(\tau,x)   -p(\tau,x,\partial)  \big) 
 v\ dx 
 \cr
&= -\int_{\partial{\caQ_\delta}}
  \Phi(\tau,x)
   \Big((V+\beta(\tau,x))u\cdot v- u\cdot (V
   + \beta(\tau,x)) v\Big)
\,d\Sigma.
\end{align}
To prove \eqref{eq:dotGreen}, subtract \eqref{eq:Aidentity}  from the same identity with 
$u$ and $v$ interchanged.

$\bullet\ $    Equations for the annihilators.
The function
$$
(\uu, \ug_1,\ug_2)\ \in\
 C^\infty(\overline{{\caQ_\delta}})
 \times
 C^\infty(\partial{\caQ_\delta};\caE^-)
 \times
  C^\infty(\partial{\caQ_\delta};\caE^+)
 $$ 
 annihilates the range if and only if $\forall u\in H^2({\caQ_\delta})$,
\begin{equation}
\label{eq:annih}
\begin{aligned}
\int_{{\caQ_\delta}} \big(  
 \tau^2\Pi    (\tau,x)   -   p  &  (\tau,x,\partial)  \big) u  \cdot \uu\ dx 
 \ +\
 \int_{\partial{\caQ_\delta}}
\pi^-(\widetilde \nu)   u \cdot \ug_1
 \ d\Sigma
  \cr
  \ &+\, 
  \int_{\partial{\caQ_\delta}}
\pi^+(\widetilde \nu)
\Big(
V+ \tau+2H_\uOmega
\Big)u
 \cdot \ug_2
 \, d\Sigma
 \  =\ 0.
 \end{aligned}
 \end{equation}
 
 The operator $\tau^2 \Pi(\tau,x) -p$ is equal to its own transpose.  Therefore,
taking $u$ that vanish on a neighborhood of $\partial{\caQ_\delta}$ implies that
\begin{equation}
\label{eq:uupde}
  \big(  
 \tau^2\Pi(\tau,x)   -p(\tau,x,\partial)  \big) \uu\ =\ 
 0
 \qquad
 {\rm on}
 \quad
 {\caQ_\delta}\,.
 \end{equation}
 This together with 
 \eqref{eq:dotGreen}
 shows that
 \eqref{eq:annih} holds if and only if 
  \begin{equation}
 \label{eq:annih2}
 \begin{aligned}
0
\ =\   &\int_{\partial{\caQ_\delta}}
  \pi^+(\widetilde \nu)
\Big(
V+ \tau+2H_\uOmega
\Big)u
\, \cdot\, \ug_2
 \,
   +\,
\pi^-(\widetilde \nu)u \cdot \ug_1
  \cr
&  -
\Phi(\tau,x)
   \big((V+\tau +2H_{\uOmega})u\cdot \uu- u\cdot (V
   + \tau + 2H_{\uOmega})\uu\big)
 \  d\Sigma.
 \end{aligned}
 \end{equation}

 Equation \eqref{eq:annih2} is used first on test functions
 $u$ that satisfy $(V+\tau + 2H_{\uOmega} )u=0$  on 
 $\partial {\caQ_\delta}$.   That constraint leaves $u|_{\partial{\caQ_\delta}}$ arbitrary.
 Of those test functions first consider those that satisfy  $\pi^-(\widetilde\nu)u|_{\partial{\caQ_\delta}}=0$.  For those
 one finds
 \begin{align*}
   \int_{\partial{\caQ_\delta}}
 \Phi(\tau,x) 
   \,
    u\cdot (V
   + \tau + 2H_{\uOmega})\uu
\ d\Sigma
 \ =\ 0\,.
 \end{align*}
 Since the  $\Phi$ factor is scalar and nowhere vanishing 
 it follows that for arbitrary $\phi\in C^\infty(\partial{\caQ_\delta})$,
 \begin{equation*}
  \int_{\partial{\caQ_\delta}}
  \pi^+(\widetilde \nu)    \phi 
    \,\cdot\, (V
   + \tau + 2H_{\uOmega})\uu
\
d\Sigma
 \ =\ 0 \,.
 \end{equation*}
 This shows that $\uu$ satisfies the transposed boundary condition
 \begin{equation}
 \label{eq:uuNeum}
 \pi^+(\widetilde\nu)^\dagger  \, 
\big(V
   + \tau + 2H_{\uOmega} \big) \uu
   \ =\ 0,
   \qquad{\rm on}\quad
   \partial{\caQ_\delta}.
   \end{equation}
   
Next take $u$ satisfying $ \pi^+(\widetilde \nu)    u     |_{\partial{\caQ_\delta}} =0$.
Then
$u|_{\partial{\caQ_\delta}} = \pi^-(\widetilde \nu)u$.  This yields
$$
  \int_{\partial{\caQ_\delta}}
  \Phi(\tau,x)
   \Big(
   \pi^-(\widetilde \nu)u\cdot (V
   + \tau + 2H_{\uOmega})\uu\Big)
   \ +\
   \pi^-(\widetilde\nu)u \cdot \ug_1
\,d\Sigma\,.
$$
The set of  functions $\pi^-(\widetilde \nu)u|_{\partial{\caQ_\delta}}$ includes
the set of 
  $\pi^-(\widetilde\nu)\psi$ for an arbitrary
$\psi\in C^\infty(\partial{\caQ_\delta};\CC^2)$.  It follows that on $\partial{\caQ_\delta}$,
\begin{equation}
\label{eq:2nd}
\pi^- (\widetilde\nu)^\dagger
\bigg(
 \Phi(\tau,x)
    \Big(V
   + \tau + 2H_{\uOmega}\Big)\, \uu
   \ +\ \ug_1
\bigg)=0
\quad
{\rm on}
\quad
\partial{\caQ_\delta}.
\end{equation}
   
Next extract the information from test functions
that satisfy $u|_{\partial{\caQ_\delta}}=0$.  For such
test functions,  $\big[ V + \tau + 2H]_{\partial{\caQ_\delta}}$
can be chosen as an arbitrary element $\psi\in C^\infty(\partial{\caQ_\delta};\CC^2)$.
This yields
 \begin{equation*}
 \label{eq:annih3}
 \begin{aligned}
  -\int_{\partial{\caQ_\delta}}
  \Phi(\tau,x)
  \,
   \psi\cdot \uu
   \ d\Sigma
   +
  \int_{\partial{\caQ_\delta}}
\pi^+(\widetilde \nu)
\psi 
\, \cdot\, \ug_2
 \, d\Sigma
 \ =\
 0.
 \end{aligned}
 \end{equation*}

First take those $\psi$ that satisfy 
$\pi^+(\widetilde\nu)\psi=0$.  That is equivalent to $\psi = \pi^-(\widetilde\nu)\phi$
for arbitrary $\phi$.  That yields
$$
\int_{\partial{\caQ_\delta}}
  \Phi(\tau,x)
  \,
   \pi^-(\widetilde\nu)\phi\cdot \uu
   \ d\Sigma
\ =\ 0\,.
$$
This is equivalent to the Dirichlet boundary condition for $\uu$,
\begin{equation}
\label{eq:uuDir}
\pi^-(\widetilde \nu)^\dagger\uu
\ =\ 0,
\quad
{\rm on}
\quad
\partial{\caQ_\delta}.
\end{equation}

Finally, consider $\psi$ with $\pi^-(\widetilde\nu)\psi=0$.  Equivalently $\psi = \pi^+(\widetilde\nu)\phi$
for arbitrary $\phi$.  This yields
$$
 \int_{\partial{\caQ_\delta}}
\pi^+(\widetilde\nu)\phi \cdot
\Big(
- \Phi(\tau,x)
  \uu
\ +\ \ug_2
\Big)
 \ d\Sigma
 \  =\ 0.
 $$
Since $\phi$ is arbitrary this is equivalent to 
\begin{equation}
\label{eq:fourth}
\pi^+
(\widetilde \nu)^\dagger\Big(
- \Phi(\tau,x)
\uu
\ +\  \ug_2\Big) = 0
\quad
{\rm on}
\quad
\partial{\caQ_\delta}.
\end{equation}

$\bullet$  Proof that $\uu=0$, $\ug_1=0$, and $\ug_2=0$.  The three equations \eqref{eq:uupde}, \eqref{eq:uuNeum},  and
\eqref{eq:uuDir}
 assert that $\uu$ is a smooth solution of
the transposed  boundary value problem with zero sources.  Theorem \ref{thm:perturbed}
implies that $\uu=0$.

From the fact that $\uu =0$, \eqref{eq:2nd} implies that
$(\pi^-(\widetilde \nu))^\dagger \ug_1=0$.  In addition $\ug_1$ takes values in $\caE^-(\widetilde \nu)$.
There is an  $R_2$ so that for $|\tau|>R_1$,  
$\pi^-(\widetilde\nu)^\dagger$ is injective on
$\caE^-(\widetilde\nu)$ for all $x\in \partial {\caQ_\delta}$.  
For those $\tau$, conclude that $\ug_1=0$.

An entirely analogous argument using \eqref{eq:fourth} shows
that $\ug_2=0$.  This completes the proof that the
annihilator of the range is equal to $\{0\}$.
\hfill
\qed
\vskip.2cm

\subsection{Analyticity in $\tau$ of the Helmholtz solution}

Use the shorthand $\caE^\pm(\tau,x)$ for 
$\caE^\pm(\widetilde\nu(\tau,x) ) $.
The vector spaces
$\caE^\pm(\tau,x)$ depends analytically on $\tau$.
The next example shows that 
defining what it means   to depend
analytically on $\tau$ has pitfalls.

\begin{example}  
\label{ex:holsub}
{\bf i.}  The subspace $\UU(\tau)\subset \CC^2$ spanned by $(1,\tau^2)$ 
  depends analytically on $\tau$
for any reasonable definition including the one  below.

{\bf ii.}   The unit vectors spanning
$\UU(\tau)$ 
are 
$$
e^{i\theta(\tau)}\ 
\frac{(1,\tau^2)}
{(1+|\tau|^4)^{1/2}},
\qquad
\theta\in \RR.
$$
No choice of $\theta$ makes this holomorphic.

{\bf iii.}  Orthogonal projection onto $\UU(\tau)$ has matrix 
equal to
$$
\frac1{1+|\tau|^4}
\begin{pmatrix}
1 & \overline \tau^2
\cr
\tau^2 & |\tau|^4
\end{pmatrix} .
$$
It is not a holomorphic function of $\tau$.
\end{example}

The analytic dependence of $\caE^\pm(\tau,x)$ 
is expressed as follows.  
 For each $(\tau,x)$, $\CC^2=\caE^+(\tau,x) \oplus \caE^-(\tau,x)$.
For $\tau$ near a fixed $\utau$ and all $x\in \partial{\caQ_\delta}$, $\pi^+(\widetilde\nu)$ is an isomorphism from
$\caE^+(\tau,x)\to \caE^+(\utau,x)$.  Define the linear transformation $R^+(\tau,x)\in {\rm Hom}(\CC^2)$
to be the inverse of this isomorphism for $v\in \caE^+(\utau,x)$ and
equal to zero on $\caE^-(\utau,x)$.  An analogous definition yields
$R^-(\tau,x)$.    Then $R^\pm(\tau,x)\in {\rm Hom}(\caE^\pm(\utau,x):\CC^2)$ depend analytically on $\tau$.
For $\tau$ near $\utau$ and all $x\in \partial{\caQ_\delta}$,
$$
\caE^+(\tau,x) \ =\ R^+(\tau,x)\, \caE^+(\utau,x)\,.
$$
\JR{
For any $\ux\in \partial{\caQ_\delta}$ one can choose a nonzero element $\bfe\in \caE^+(\utau,\ux)$.
Then for $x$ in a neighborhood of $\ux$ and $\tau$
in a neighborhood of $\utau$, $R^+(t,x)\bfe$ is a smooth basis of 
$\caE^+(\tau,x)$  that depends holomorphically on $\tau$.   The existence of 
such a local basis is what it means for the $\tau$ dependent vector bundles $\caE^+(\tau,x)$
to be holomorphic in $\tau$.
The problem with the choice in  {\bf ii,iii} of Example \ref{ex:holsub} is caused by the 
normalizations.
}

The boundary value problem \eqref{eq:PHBVP}  
has source terms $g_j$ that take values in 
$\caE^\pm(\tau,x)$.   
The local representation allows one to suppress the 
$\tau$ dependence as follows.  For $\tau$ near $\utau$,
a section $g$ of $\caE^+(\tau,x)$ is uniquely represented as
$R^+(\tau,x)\ug$ where $\ug$ is takes values in the $\tau$ dependent
space
$\caE^+(\utau,x)$.  The boundary value problem takes the form
\begin{equation}
\label{eq:PHBVP2}
\begin{aligned}
\big(
\tau^2\,\Pi(\tau,x)\ -\ 
p(\tau,x,\partial) 
\big)  u \ &=\ f\,
\quad 
{\rm on}
\quad
{\caQ_\delta},
\cr 
\pi^+(\widetilde \nu(\tau,x))  u &= R^-(\tau,x)\ug_1
\quad
{\rm on}
\quad
\partial{\caQ_\delta},
\cr
\pi^+(\widetilde \nu(\tau,x))
\big(
V+ \tau+2H_\uOmega(X(\tau,\ux))
\big)u
 &=   R^+(\tau,x)\ug_2
\ \
{\rm on}
\ \
\partial{\caQ_\delta}.
\end{aligned}
\end{equation}
Here  $\ug_1$ takes values in $\caE^-(\utau,x)$ and 
$\ug_2$ takes values in $\caE^+(\utau,x)$.
In this form, the source terms $\ug_j$ belong to a $\tau$-independent
space and the coefficients of the operators depend
differentiably on $\tau,x$ and analytically on $\tau$.   

\begin{definition}
A $\tau$-dependent section $g_1(\tau)\in H^{3/2}(\caE^-(\tau,x))$ 
depends analytically on $\tau$ when the corresponding functions
$\ug_1(\tau) \in H^{3/2}(\caE^-(\utau,x))$ depend analytically
on $\tau$.
A similar definition applies for $g_2(\tau)\in H^{1/2}(\caE^-(\tau,x))$.
\end{definition}

\begin{theorem} 
\label{thm:holomorphy}  If the source terms 
$$
(f,g_1,g_2)
\ \in\
  L^2({\caQ_\delta}) \times H^{3/2}(\caE^-(\tau,x)) \times H^{1/2}(\caE^+(\tau,x)) 
  $$
   depend analytically
on $\tau$ on ${\rm Re}\,\tau>M$, then the corresponding solution
$u(\tau,\cdot)$
of \eqref{eq:PHBVP}  is an analytic function of $\tau$ with values in 
$H^2({\caQ_\delta})$.
\end{theorem}

{\bf Proof.}  Standard elliptic theory shows that writing
$\tau=a+ib$  the map $a,b\mapsto u$ is infinitely differentiable
with values in $H^2({\caQ_\delta})$.
The derivatives 
satisfy the system obtained by differentiating, with respect to $a,b$,
the system and boundary conditions satisfied by $u$.   
 
To prove analyticity it suffices to show that  $w:= \partial u/\partial \overline \tau=0$.
Since all the coefficients and the  $f,g_1, g_2$ are analytic,
differentiating the boundary value problem with respect to $\overline \tau$
shows that $w$ satisfies
\begin{equation*}
\begin{aligned}
\Big(
\tau^2\,\Pi(\tau,x)\ -\ 
p(\tau,x,\partial) 
\Big)  w \ &=\ 0\,
\quad 
{\rm on}
\quad
\caQ_\delta,
\cr 
\pi^-(\widetilde \nu(\tau,x)) w \ &=\ 0
\quad
{\rm on}
\quad
\partial\caQ_\delta,
\cr
\pi^+(\widetilde \nu(\tau,x))
\Big(
V(\tau,x,\partial)+ \tau+2H_\uOmega(X(\tau,\ux))
\Big)w
\ &= \ 0
\quad
{\rm on}
\quad
\partial\caQ_\delta.
\end{aligned}
\end{equation*}
Theorem \ref{thm:perturbed},
 implies that $w=0$.
\hfill
\qed
\vskip.2cm

\section{Proofs of the Main Theorems}
\label{sec:maintheorems}

\JR{The main elements of the proofs of the Main Theorems have been
prepared.   In this section they are combined to finish the proofs.}

\subsection{The stretched equation on $\caQ_\delta$, \JR{proof of} Theorem \ref{thm:deltastretched}}
\label{sec:stretcheddelta}

{\bf Proof of Theorem \ref{thm:deltastretched}.}  
{\bf  Uniqueness.}  
Multiply the differential equation $L(\tau,\widetilde\partial)u^\delta =F$ 
from 
\eqref{eq:stretchedpair}
by $\Pi(\tau,x)\big(\tau-\sum A_j\widetilde\partial_j)$ and 
use
\eqref{eq:prodpauli}
to find the first  line in the Helmholtz boundary value problem
\begin{equation}
\begin{aligned}
\label{eq:helm1}
\big(\tau^2\Pi(\tau,x) \ -\ p(\tau,x,\partial)\big)u^\delta\ &=\
\Pi(\tau,x)\big(\tau-\sum A_j\widetilde\partial_j)F\,,
\cr
\pi^-(\widetilde\nu)\,
u^\delta\ &=\ 0,
\quad
{\rm on}\quad 
\partial\caQ_\delta,
\cr
\pi^+(\widetilde \nu)
\big(
V(\tau,x,\partial) + \tau + 2H_{\underline \caQ_\delta} \big)
u^\delta\ &=\ 0,
\quad
{\rm on}
\quad
\partial\caQ_\delta.
\end{aligned}
\end{equation}
The second line is part of \eqref{eq:prodpauli}.
The last line follows from 
part {\bf iv} of 
Proposition  \ref{prop:analyticBC} 
since $F=0$ on a neighborhood of $\partial\caQ_\delta$ and $u^\delta\in H^2(\caQ_\delta)$.

The hypotheses of Theorem
\ref{thm:perturbed} are satisfied. Apply the estimate of that Theorem with $f=0$ to conclude that 
$u=0$.

\vskip.1cm
 {\bf  Existence.}  
For ${\rm Re}\,\tau>M$,
Theorem 
\ref{thm:taufixed} implies that the boundary value
problem 
\eqref{eq:helm1}
has a unique solution  $u^\delta\in H^2(\caQ_\delta)$.
Theorem 
\ref{thm:holomorphy}  implies that $u$ is holomorphic with values in 
$H^2(\caQ_\delta)$.

Since $F\in L^2_{\ell\overline\Q}(\RR^3)$ it follows that 
the source term $f:=\Pi(\tau,x)\big(\tau-\sum A_j\widetilde\partial_j)F$ belongs to 
$L^2(\RR^3)$ with ${\rm supp}\,F\subset \ell\overline\Q$.  Estimate
\begin{align*}
\Big|
\int_{\caQ_\delta} f \ \ou\ dx\Big| 
\ &=\ 
\Big|
\int_{\caQ_\delta} F\ \big(\tau-\sum A_j\widetilde\partial_j)^*(\Pi \ou)\ dx\Big|
\cr
\ &\lesssim\
\|F\|_{L^2_{\ell\overline\Q}(\Q_\delta)}\,
\Big(
\|\tau u \|_{L^2(\Q_\delta)}\ +\ \| \nabla_xu \|_{L^2(\Q_\delta)}\Big)\,.
\end{align*}
Estimate  the two terms on the right as follows.    Write 
\begin{align*}
&C\| \mu^{-1}  F\|_{L^2_{\ell\overline\Q}(\Q_\delta)}\ 
\|\mu \tau u \|_{L^2(\Q_\delta)}
\, \le\,
\frac
{C^2\mu^{-2} }2  \|F\|_{L^2_{\ell\overline\Q}(\Q_\delta)}^2
 +
\frac{\mu^2|\tau|^2}2
\|u \|_{L^2(\Q_\delta)}\big).
\cr
&C\|\eps^{-1}F\|_{L^2_{\ell\overline\Q}(\Q_\delta)}\,
 \|\eps \nabla_xu \|_{L^2(\Q_\delta)}
\ \le\
\frac{C^2\eps^{-2}}2  \|F\|_{L^2_{\ell\overline\Q}(\Q_\delta)}^2
 +
\frac{\eps^2}2 |\|\nabla_xu \|_{L^2(\Q_\delta)}^2.
\end{align*}
Choose $\mu,\eps$ so that
$\mu^2|\tau|^2=|\tau|(\Re\tau)$ 
and
$\eps^2 = (\Re\tau)/ |\tau|$.
Then,
\begin{align*}
&C\|   F\|_{L^2_{\ell\overline\Q}(\Q_\delta)}\ 
\| \tau u \|_{L^2(\Q_\delta)}
\ \le\
\frac{C^2|\tau|}{2\Re\tau}\,
  \|F\|_{L^2_{\ell\overline\Q}(\Q_\delta)}^2
 +
\frac{|\tau|(\Re\tau)}{2}\, \|u \|_{L^2(\Q_\delta)}^2.
\cr
&C\|F\|_{L^2_{\ell\overline\Q}(\Q_\delta)}\,
 \| \nabla_xu \|_{L^2(\Q_\delta)}
\ \le\
\frac{C^2|\tau|}{2\Re\tau}  \|F\|_{L^2_{\ell\overline\Q}(\Q_\delta)}^2
 +
\frac{\Re\tau}{2|\tau|}\, \|\nabla_xu \|_{L^2(\Q_\delta)}^2.
\end{align*}

Absorbing the two right hand terms, 
 Theorem \ref{thm:perturbed}
 shows that with constant independent of $\delta$,
\begin{equation}
\begin{aligned}
\label{eq:h1est2}
|\tau|\,({\rm Re}\, \tau)\, \big\|{u^\delta} &  \big\|_{L^2(\caQ_\delta)}^2 
\ +\
|\tau|\, 
\big\|{u^\delta}\big\|_{L^2(\partial\caQ_\delta)}^2
\cr
& +\
\frac{\Re\tau}{|\tau|}\,\big\|\nabla {u^\delta}\big\|_{L^2(\caQ_\delta)}^2\ \lesssim\
 \frac{|\tau|}{\Re\tau}\,
\big\| {F}\big\|_{L^2_{\ell\overline\Q}(\caQ_\delta)}^2.
\end{aligned}
\end{equation}
Multiplying by $(\Re\tau)/|\tau| $ yields 
\eqref{eq:elliot}.

To complete the proof it suffices to show that $u^\delta$ satisfies the 
stretched boundary value problem
 \eqref{eq:bermagic10} on $\caQ_\delta$.     To
do that,  reverse the steps that lead from the stretched equations to the
Helmholtz boundary value problem.
The proof uses the $H^2(\Q_\delta)$ regularity that requires the smoothness of
$\Q_\delta$.  Define
$$
w\ :=\ 
\big(
A(\widetilde \partial)
\,+\, \tau
\big)u\ \in \ H^1(\caQ_\delta).
$$
It suffices  to show that the stretched equation, $w=F$, is satisfied.

The Helmholtz equation implies that 
$w-F\in H^1(\caQ_\delta)$ satisfies 
\begin{equation}
\label{eq:penny1}
\big(
A(\widetilde \partial)
\,-\, \tau
\big)
\Big(w\ -\ F\Big)
\ =\
0\,,
\qquad
{\rm on}
\quad
\caQ_\delta\,.
\end{equation}

Part {\bf iv} of Proposition  \ref{prop:analyticBC}
 shows that the  derivative boundary condition satisfied by 
 $u^\delta$  is equivalent to 
$$
\pi^+\big(A(\widetilde\nu)\big)\ 
A(\widetilde\nu)^{-1}\
\big(
w-F
\big)
\ =\ 0
\qquad
{\rm on }
\quad
\partial
\caQ_\delta\,.
$$
Since $\pi^+(A(\widetilde \nu))$ and $A(\widetilde \nu)$ commute,
this is equivalent to 
\begin{equation}
\label{eq:penny2}
\pi^+\big(A(\widetilde\nu)\big)\ 
\big(
w-F
\big)
\ =\ 0
\qquad
{\rm on }
\quad
\partial
\caQ_\delta\,.
\end{equation}

When $\tau$ is  real and large, 
the pair of equations 
\eqref{eq:penny1}, \eqref{eq:penny2}
is a strictly dissipative
boundary value problem with vanishing sources on the smooth
domain $\caQ_\delta$
with noncharacteristic boundary.  
The solution is in $H^1(\Q_\delta)$.    That the solution vanishes follows 
by a direct integration by parts showing that 
$$
\|w-F\|_{L^2(\Q_\delta)}^2
 \ \lesssim\
 {\Re}\,
 \int_{\Q_\delta} \Big(
 \big(\tau- A(\partial)\big)
 (w-F), (w-F)
 \Big)_{\CC^6}\ dx
 \ =\ 0\,.
 $$
 
The map $\tau\mapsto (w-F)(\tau)$
 is 
holomorphic for ${\rm Re}\,\tau $ large.
It has just been proved that 
it vanishes on $]m,\infty[$ for $m$ large.
By analytic continuation,
 it follows that $w-F=0$ for all ${\rm Re}\,\tau >M$.

 Thus the stretched equation is satisfied on $\caQ_\delta$
 for ${\rm Re}\,\tau>M$.
This  completes the proof of existence.
\hfill
\qed
\vskip.2cm

\subsection{The stretched equation on $\caQ,$
\JR{proof of}
 Theorem \ref{thm:mainstretched}}
\label{sec:stretched}

{\bf Proof of Theorem \ref{thm:mainstretched}.}
{\bf  Uniqueness.}   
The  solution with vanishing data
is holomorphic in 
${\rm Re}\,\tau$ large.  To prove that it vanishes it is sufficient to prove that
it vanishes 
for $\tau\in ]m,\infty[$ for $m$ large.

  For $\tau$ real and large, the stretched equation,
 $L\big(\tau,\widetilde \partial\,\big)    u(\tau)=0$
is symmetric positive in the sense of Friedrichs, that is
$$
L\big(\tau,\widetilde \partial\,\big)  + L\big(\tau,\widetilde \partial\,\big)^* 
\ \ge\ 
C_1( \tau -  C_2)I,
\qquad
C_1>0\,.
$$
 In addition, $ u(\tau)\in H^1(\Q_\delta)$ satisfies  strictly dissipative
boundary conditions  on each smooth faces $G_j$
 Therefore  a straightforward integration by parts
shows that 
$$
\|  u(\tau)\|_{L^2(\Q)}^2 
\ \le\
C_1\big( \tau -  C_2\big) 
\int_{\Q}
\big(
L(\tau,\widetilde \partial)   u\,,\,   u
\big)
\ dx\ =\ 0.
$$

\vskip.1cm
{\bf Existence.}   Use Theorem \ref{thm:deltastretched}.
Solve on $\caQ_\delta$  and
pass to the limit $\delta\to 0$.  At the same time
one must smooth the source term $f$ in order to apply
Theorem \ref{thm:deltastretched}.

Choose $0< \underline\eps < {\rm dist}\,({\ell\overline\Q},\partial\caQ)/2$.   
Define $K^\prime$ to be the set of points at distance
$\underline \eps$ from $\ell\Q$.    Then $K^\prime\subset \caQ$
is compact.    
For $\eps<\underline \eps$,
define $F_\eps:= j_\eps*F$ where $j_\eps$ is a
smooth mollification kernel on $\RR^3$  with support in the ball of
radius $\eps$ at the origin.    The  source
term $F_\eps\in C^\infty_{K^\prime}(\caQ)$.
   For $\delta$ sufficiently small 
$K^\prime\subset \caQ_\delta$ and 
Theorem \ref{thm:deltastretched} applies.

Define $\delta(n)=2^{-n}$, 
and 
$u^{\delta(n)}\in H^2(\caQ_{\delta(n)})$
to be 
the solution from 
Theorem \ref{thm:deltastretched}
with source term equal to $F_{\delta(n)}$.
Then with $C$
independent of $n$,
\begin{equation}
\label{eq:un}
\begin{aligned}
   ({\Re}\, \tau)^2\, \big\|   & u^{\delta(n)}\big\|^2_{L^2(\caQ_{\delta(n) } ) }
  +
(\Re\tau)\,   
\big\|  u^{\delta(n)}   \big\|^2_{L^2(\partial\caQ_{\delta(n)} ) }  
\cr 
 &+
 \frac{(\Re\tau)^2}{|\tau|^2}
  \big\| \nabla_x u^{\delta(n)}   \big\|^2_{L^2(\caQ_{\delta(n)} ) }   
 \ \le \
 C\,
 \big\| 
 F_{\delta(n)}
\big\|_{L^2_{\ell\overline\Q}(\caQ_{\delta(n)}) }^2
 \,.
\end{aligned}
\end{equation}

Extract a subsequence that converges weakly in 
$H^1(\caQ_{\delta(1)})$ to a limit $v_1$.
Extract a further subsequence that 
converges weakly in 
$H^1(\caQ_{\delta(2)})$ to a limit $v_2$.
And so  forth.  
For each $n>1$, one has $v_n= v_{n-1}$ on $\caQ_{\delta(n-1)}$.
Define $v\in H^1(\caQ)$ by $v=v_n$ on $\caQ_{\delta(n)}$.   
Using that
$\caQ_{\delta(n)}\nearrow\caQ$ and  $\partial\caQ_{\delta(n)}\cap G_j \nearrow G_j$
conclude that 
for each $n$,  $u_k$ converges 
weakly to $v$ in $H^1(\caQ_{\delta(n)})$ with
\begin{equation}
\label{eq:vO}
\begin{aligned}
 ({\Re}\, \tau)^2\, \big\|v\big\|^2_{L^2(\caQ ) }
 \ +\ &
(\Re\tau)\,   
\big\|  v  \big\|^2_{L^2(\partial\caQ ) }  \cr
\ &+\ 
\frac{(\Re\tau)^2}{|\tau|^2}\,  \big\| \nabla_x v  \big\|^2_{L^2(\caQ ) }   
 \le 
C\big\| F\|_{L^2_{\ell\overline\Q}(\caQ)}^2\,.
\end{aligned}
\end{equation}

The differential equation $L(\tau,\widetilde\partial)v =F$
on $\caQ$ follows from the  equations
$L(\tau, \widetilde\partial) u_k = F_k$ on $\caQ_{\delta(n(k))}$
on passing to the limit $k\to \infty$.   Similarly, the boundary condition
$$
\pi^+(\nu) v \ =\ 0, \qquad
{\rm on}
\quad
G_k
$$
follows on passing to the limit in 
$$
\pi^+(\nu)\,u_{\delta(n)}\big|_{G_k\cap\partial\caQ_{\delta(n)}  } 
\ =\ 0.
$$

For any $\udelta>0$  the holomorphy of $\tau \mapsto v(\tau)$ from ${\rm Re}\,\tau>M$
to $L^2(\caQ_\udelta)$ 
  follows from the fact that it is
 the weak limit of bounded family of 
holomorphic functions.   Therefore, for any  $\delta$,
$v:\{{\rm Re}\,\tau>M\}\to L^2(\caQ_\delta)$ is holomorphic.

To show that $v$ is holomorphic with values in $L^2(\caQ)$
it is sufficient to show that $\tau\mapsto \ell(v(\tau))$
is holomorphic for each $\ell$ in the dual of of $L^2(\caQ)$.

Since 
$v \in 
L^\infty\big(
\{ {\rm Re}\,\tau >M\}
\,;\,L^2(\caQ)\big)$,
it
 suffices to show that $\ell(v(\tau))$ is holomorphic 
for $\ell$ in a dense subset.    Indeed if $\ell$ is the limit
of $\ell_j$ for which the result is true, estimate
$$
\big|
\ell(v(\tau)) \ -\ \ell_j(v(\tau))
\big|
\ \le\ 
\| \ell-\ell_j\|
\ 
\sup_{{\rm Re}\,\tau>M}\ 
\| v(\tau)\|_{H^1(\caQ)},
\quad
{\rm on}\ \ 
{\rm Re}\,\tau>M.
$$
This proves that $\ell(v(\tau))$ is the uniform limit of the 
holomorphic functions $\ell_j(v(\tau))$.

Take the dense set to be the linear functionals
$v\mapsto \int v\cdot \phi\,dx$ with $\phi\in C^\infty_0(\caQ)$.
For each such $\phi$,
$\phi\in C^\infty_0(\caQ_\delta)$ for $\delta$ small.   That
$\ell(v(\tau))$ is holomorphic 
then follows from the fact that
$v$ is holomorphic with values in $L^2(\caQ_\delta )$.
This completes the proof of the Theorem.
\hfill
\qed
\vskip.2cm

\subsection{B\'erenger's equation on $\RR_t\times \caQ$,  \JR{proof of}
Theorem
\ref{thm:pauli}    }
\label{sec:splitpauli}

\JR{
Paley-Wiener Theorem for functions with values in a 
Hilbert space $H$  (see \cite{hille1996}) is needed.}

\begin{theorem} {\bf (Paley-Wiener)}
\label{thm:laplace}
  The Laplace transforms of functions
$F\in e^{Mt}\, L^2(\RR\,;\,H)$ 
with ${\rm supp}\, F\subset\{t\ge 0\}$  
 are exactly
the functions $G(\tau)$ holomorphic in ${\rm Re}\,\tau>M$
with values in $H$ and so that
$$
\sup_{\lambda>M}
\
 \int_{ {\rm Re} \,\tau  =\lambda} \big\|\widehat F(\tau)\big\|_H^2
\ 
|d\tau|
\ <\ \infty
\,.
$$
In this case the function $\widehat F(\tau)$  has
trace at ${\rm Re}\,\tau =M$ that satisfies
$$
\int e^{ -2M t } \ \|F(t)\|_H^2\ dt
=
\sup_{\lambda>M}
\
 \int_{ {\rm Re} \,\tau  =\lambda} \big\|\widehat F(\tau)\big\|_H^2
\ 
|d\tau|
 =
 \int_{ {\rm Re} \,\tau  =M} \big\|\widehat F(\tau)\big\|_H^2
\ 
|d\tau|
\,.
$$
\end{theorem}

{\bf Proof of Theorem \ref{thm:pauli}.}
{\bf Uniqueness.}    \JR{Next}  show that if $U^1,U^2,U^3$ is a solution with source $f=0$,
then $U^j=0$.   Denote by $\widehat U^j$ the Laplace transform that is holomorphic
in $\{{\rm Re}\,\tau >M\}$ with values in $L^2(\caQ)$.

The function $v(\tau):=\sum \widehat U^j$ 
is holomorphic with values in $L^2(\caQ)$ and satisfies the stretched equation
$$
\tau  \ v
\ +\
\sum A_j \widetilde \partial_j v \ =\ 0\,.
$$
In addition, $v|_{G_k}$ is holomorphic with values in 
$L^2(G_k)$.   The boundary  condition satisfied by $\sum U^j$ implies that $v$
  satisfies the  boundary condition
$$
v|_{G_k} \ \in \ 
\caE^+(\nu)\,,
\qquad
1\le k\le 6.
$$
The stretched operator is elliptic.  When $\tau\in ]m,\infty[$ the stretched
operator is symmetric and positive in the sense
of Friedrichs.  The 
uniqueness theorem for such strictly dissipative symmetric and elliptic problems with 
trihedral corners from  Part I of 
\cite{HR:2016:HBV}
implies that $\widehat u(\tau)=0$ for $\tau\in ]m,\infty[$.  By analytic continuation,
$v(\tau)=0$ on $\{ \Re\tau>M\}$.

The Laplace transform of the split equation yields
$$
\big(\tau+ \sigma_1(x_1)\big) \widehat U^j \ =\  - A_1\partial_1
v
\ =\
0\,.
$$
This implies that $\widehat U^j $ vanishes and therefore that 
$U^j=0$.  This completes the proof of uniqueness.

\vskip.1cm
 {\bf Existence.}   The solution $u(t,x)$ is constructed by 
finding its Laplace transform.   Denote by 
$U^1(t,x)$, $U^2(t,x)$, and, $U^3(t,x)$  the unknowns to be found.      Denote by $v(\tau,x)$ the 
function of $\tau$ that will be   the Laplace transform of 
$U^1(t,x) + U^2(t,x) + U^3(t,x) $.
Define $v(\tau,x)$ to be the solution of 
the stretched equation
\begin{equation}
\label{eq:vsplit}
\tau  \ v
\ +\
\sum_{j=1}^3 A_j \widetilde \partial_j v \ =\ F(\tau):=
\, \sum_{j=1}^3\,
 \frac{\tau \widehat f_j(\tau)}{\tau+\sigma_j(x)}
\,.
\end{equation}
constructed in Theorem 
\ref{thm:mainstretched}.  Then $v$ holomorphic in ${\rm Re}\,\tau>M$ with values
in $H^1(\caQ)$ 
and
$v|_{G_k}$ is holomorphic with values in 
$L^2(G_k)$.
In addition,
\begin{equation}
\label{eq:lapest}
\begin{aligned}
({\rm Re}\, \tau)
\big\| v(\tau)
\big\|_{L^2(\caQ)}
 \, +\, &
({\rm Re}\, \tau)^{1/2}
\big\|\,
 v(\tau)
\big\|_{L^2(\partial\caQ)}
\, +\, 
\frac{\Re\tau}{|\tau|}\,\big\| \nabla_xv(\tau)
\big\|_{L^2(\caQ)}
\cr
&\ \le\ 
C\, \big\| F(\tau)\big\|_{L^2_K(\caQ)}
\ \le\ 
C\, \big\| \widehat f(\tau)\big\|_{L^2_K(\caQ)}.
\end{aligned}
\end{equation}

Define 
$V^j$ destined to be the Laplace transforms of the $U^j$ by the analogue of 
\eqref{eq:split10},
 \begin{equation}
\label{eq:split20}
 \big(\tau+ \sigma_j(x_j)\big) V^j \ +\ A_j\partial_j
v= \ \widehat {f_j}
\,,
\qquad
j=1,2,3\,.
 \end{equation}
Multiplying  by  $\tau/(\tau+\sigma_j(x_j))$
yields
\begin{equation*}
\tau\,
V^j + A_j
\widetilde{\partial_j}v
\ =\
 \frac{\tau \widehat {f_j}}{    \tau+\sigma_j }
\,,
\quad
j=1,2,3\,.
\end{equation*}
Summing yields
\begin{equation*}
\label{eq:othersplit}
\tau  \Big( V^1 + V^2 + V^3\Big) 
\ +\
\sum A_j \widetilde \partial_j v \ =\
\
\sum_{j=1}^3
 \frac{\tau\widehat {f_j}}{ \tau+\sigma_j(x_j)   }=F\,.
\end{equation*}

 Subtracting from \eqref{eq:vsplit}  yields
\begin{equation*}
\label{eq:sumV}
\tau  \Big(V^1 + V^2 + V^3-v\Big) 
 \ =\ 0
 \quad
 {\rm so},
 \quad
 v=V^1+ V^2+ V^3.
\end{equation*}

The Paley-Wiener theorem implies that 
$$
\sup_{\lambda>M}
\int \| \widehat f(\tau)\|^2\ |d\tau|\ \le \
\int e^{2Mt} \|f(t)\|_{L^2_{\ell\overline\Q}(\caQ)}^2
\ dt.
$$
Equation \eqref{eq:lapest} together with the Paley-Wiener Theorem
 implies that  $v$ is the Laplace transform of a function
 $u\in e^{Mt} L^2(\RR;L^2(\caQ))$ supported in $t\ge 0$.
Moreover,
\begin{equation*}
\begin{aligned}
\int_0^\infty
e^{2Mt} \Big(
M \big\| u(t)
\big\|_{L^2(\caQ)}^2
 +  
M^{1/2}
\big\|\,
 u(t)|_{\partial O}   &
\big\|_{L^2(\partial\caQ)}^2
\Big)\,dt
\cr
& \lesssim 
\int_0^\infty e^{2Mt} \|f(t)\|_{L^2_{\ell\overline\Q}(\caQ)}^2
\, dt.
\end{aligned}
\end{equation*}

Similarly the Paley-Wiener Theorem implies that 
$V^j(\tau)$ is the Laplace transform of a function
$U^j(t)\in e^{Mt}L^2(\RR;H^{-1}(\caQ))$ supported
in $t\ge 0$ and satisfying
$$
\int_0^\infty
e^{2Mt} \,
\big\| MU^j(t), \partial_t U^j(t)
\big\|_{H^{-1}(\caQ)}^2
\ dt\
\lesssim \
\int_0^\infty e^{2Mt} \|f(t)\|_{L^2_{\ell\overline\Q}(\caQ)}
\, dt.
$$

The fact that $v=\sum V^j$  implies that $u=\sum U^j$.
Equation \eqref{eq:split20}
 implies  that 
$(U^1,U^2,U^3)$ 
satisfies the 
B\'erenger split equations.
The last two estimates are exactly those required in 
Theorem \ref{thm:pauli}.

Denoting by  $\LL$ the Laplace transform, one has
$$
\LL
\big(\pi^-(\nu)\,u|_{G_j} \big)
\ =\ 
\pi^-(\nu) \big(
\LL(u|_{G_j} ) \big)
\ =\
\pi^-(\nu) \, v|_{G_j}
\ =\ 0\,.
$$
This proves the boundary condition 
$\pi^-(\nu)\,u|_{G_j} =0$.
This completes the proof that the $U^j$ satisfy the 
boundary value problem and estimates of Theorem \ref{thm:pauli}.
\hfill
$\Box$

\bibliographystyle{abbrv}
\bibliography{pml}

\end{document}